%% file: 25-note-coordination-etc-stc_2column.tex
 \documentclass[twocolumns]{IEEEtran}

\usepackage{graphicx,fleqn,amstext,amsmath,amssymb,color,psfrag} 
\usepackage[latin1]{inputenc}

\usepackage[colorlinks=true, letterpaper,    
      urlcolor=black,     
      citecolor=blue,
      menucolor=black,    
      linkcolor=black,    
      pagecolor=black,    
      breaklinks=true,
]{hyperref}

\include{commands}

\begin{document}

\title{\LARGE A Lyapunov redesign of coordination algorithms for cyber-physical systems}
\author{Claudio De Persis and Romain Postoyan
\thanks{Claudio De Persis is with the Faculty of Mathematics and Natural Sciences, University of Groningen, the Netherlands, {\tt\small c.de.persis@rug.nl}.
His work is partially supported by the The Netherlands Organization for Scientific Research (NWO) under the auspices of the project QUICK (QUantized Information
Control for formation Keeping).}
\thanks{R. Postoyan is with the Universit\'e de Lorraine, CRAN, UMR 7039 and the CNRS, CRAN, UMR 7039, France {\tt\small romain.postoyan@univ-lorraine.fr}. His work is partially supported by the ANR under the grant COMPACS (ANR-13-BS03-0004-02).}
}

\maketitle

\IEEEpeerreviewmaketitle

\begin{abstract} We investigate the coordination of a network of agents in a cyber-physical environment. In particular, we consider nonlinear agents' dynamics of arbitrary dimensions, which satisfy a strict passivity property. The objective is to ensure the convergence of the differences between the agents' output variables to a prescribed compact set (hence covering rendez-vous and formation control as specific scenarios),  while taking into account the communication and/or computation limitations to which are subject the agents. We develop event-based sampling strategies for that purpose by following an emulation approach:
 we start with distributed controllers which solve the problem in continuous-time, and we then explain how to implement these using event-based sampling.
The idea is to define a triggering rule per edge using an auxiliary variable whose dynamics only depends on the local variables. The triggering laws are designed to compensate for the perturbative term introduced by the sampling, a technique that reminds of Lyapunov-based control redesign.
All strategies guarantee the existence of a uniform minimum amount of times between any two edge events. The analysis is carried out within the framework of hybrid systems and an invariance principle is used to conclude about coordination.
\end{abstract}

\section{Introduction}\label{sect-intro}
Recent years have witnessed a massive amount of work on large-scale systems that interact locally to achieve a general coordination task.
In fact many engineered systems  have large dimensions  and requiring the different components (or agents) of these large-scale systems to exchange information only with neighboring units is valuable because it improves scalability and robustness in case of faults. On the other hand,  latest technological advances  are enabling scenarios in which computing and communication devices are an integral part of the physical processes to control. Despite this, most coordination algorithms ignore  the features of these devices, while they  may severely impact the desired agreement property. It is therefore essential to develop control strategies that take these constraints into account in their design. The problem can be addressed via the construction of event-based sampling strategies, see \emph{e.g.}, \cite{Dimarogonas-et-al-tac2012,Nowzari-Cortes-2012(aut),De-Persis-Frasca-TAC'13,De-Persis-et-al-cdc2013(self),Seyboth-et-al-aut2012}.  The idea is that each agent updates its control input \emph{only} at a sequence of time instants which depends on the local variables, and not continuously. In that way, the energy expenditure of the actuators batteries is reduced, the actuators wear is slowed down, and the usage of the computation and/or communication resources can be limited, according to the type of implementation.

Several event-based sampling paradigms exist in the literature depending on the way the sequence of input updates is defined: \emph{event-triggered control} (\cite{Arzen-99,Astrom-Bernhardsson-IFAC99}), \emph{self-triggered control} (\cite{Velasco-Fuertes-Mari-03}),  \emph{time-triggered control} (see Section \ref{sect-pb-statement} for a more detailed discussion).
These paradigms have been first proposed for single systems with a single feedback loop (see the survey \cite{Heemels-Johansson-Tabuada-cdc12} and the references therein). The multi-agent systems, on the other hand, are particularly challenging in this context.

First, these systems are generally \emph{distributed} as each agent has only access to its own state and the state of its neighbours (and not to the state of the overall system). Hence, it is necessary to design distributed triggering conditions which only depend on the local variables. One of the main difficulties here is to to ensure the existence of a minimum strictly positive amount of time between two successive triggering instants. The existence of such a time is essential for the controller to be realizable, as the hardware cannot tolerate arbitrarily close-in-time updates, as well as to rule out Zeno phenomenon.
Second, the stability analysis often relies on a weak Lyapunov function, in the sense that the derivative of the Lyapunov function along the system solution is non-positive (as opposed to strong Lyapunov functions for which it is strictly negative -- outside the attractor). This is an important difference with the vast majority of centralized stabilizing event-triggered control techniques, which require the knowledge of a strong Lyapunov function.  This point induces non-trivial technical difficulties, which also makes existing centralized event-triggering results not trivially applicable for multi-agent systems.

Despite these difficulties, several event-based algorithms have been presented for the synchronization of multi-agent systems, considering event- and self-triggered control strategies (see  \cite{OD:JL:ADHS12,Dimarogonas-et-al-tac2012,Fan-et-al-aut13,garcia-et-al-ijc13,Liuzza-et-al-nolcos2013,Nowzari-Cortes-2012(aut),De-Persis-Frasca-TAC'13,Seyboth-et-al-aut2012} to cite a few). The number of works on the topic has been growing exponentially since the appearance of  \cite{Dimarogonas-et-al-tac2012} and we do not aim at including an exhaustive survey of all the contributions.  Nonetheless, it has to be noted that most results concentrate on specific agents' dynamics, typically single- or double-integrators. The work in \cite{Liuzza-et-al-nolcos2013} is one of the rare studies which deal with agents modeled by nonlinear systems: it addresses a particular type of interconnected feedback linearizable systems. We see that there  is currently a gap between the existing techniques for the coordination of nonlinear systems in continuous-time, and their implementation in a cyber-physical environment.


In this paper, we consider a network of strictly passive systems which can have nonlinear dynamics and be of arbitrary dimensions. 
 Note  that passivity takes an outstanding role in problems of coordination control (see \emph{e.g.}, \cite{Bai2011,Burger2013,burger.depersis.aut13,depersis.jaya.TCNS14,schaft.maschke.sicon13}).
Our objective is to design distributed controllers which ensure that the difference between the agents' outputs -- which we call \emph{relative distances} -- converge to a prescribed compact set, as in \cite{Arcak-tac07}. This general formulation encompasses rendez-vous and formation control as particular cases, and can be extended to deal with several cooperative control problems. To our purpose, we follow an emulation approach as
we start from the distributed controllers proposed in \cite{Arcak-tac07}, which solve the problem in continuous-time, and we then design a triggering condition per edge to decide when to update the corresponding control input. To do so, we start from an energy-like Lyapunov function from \cite{Arcak-tac07} 
and we add a term that takes into account the `energy' associated with the sampling error. This addition is necessary to overcome the occurrence of extra terms that would disrupt the convergence of the algorithm.
We let this extra term depend on \emph{clock} variables (one per each edge in the network), which we introduce to regulate the sampling. We then synthesize the clock dynamics
in such a way that the overall Lyapunov function computed along the trajectories of the system remains monotonically decreasing despite the sampling.
We stress that, although the vast majority of the results available in event-based control of multi-agent systems is based on Lyapunov analysis and design, to the best of our knowledge this is the first time in the context of event-based control of network systems that the candidate `physical' Lyapunov function is extended to take into account the `cyber part' of the system and give rise to the triggering rule.
The idea to introduce clocks to define the triggering rule is inspired by the work on sampled-data systems in \cite{Carnevale_Teel_Nesic_netw_TAC07}, which has been adapted to event-triggered control in \cite{Postoyan-Tabuada-Nesic-Anta-CDC'11}.

We first assume that the relative distances are continuously available, in which case we derive event-triggered control laws. Afterwards, we explain how to derive (aperiodic) time-triggering rules. 
It has to be noted that these results apply to heterogonous networks (\ie the agents are not required to have the same dynamics), which is also a novelty. We then focus on homogenous networks and we develop self-triggered controllers, under an additional assumption. The existence of a uniform strictly positive lower bound on the inter-edge events is guaranteed in all cases. The overall systems are modelled as hybrid systems using the formalism of \cite{Goebel-Sanfelice-Teel-book} and the analysis  invokes an invariance principle from \cite{Goebel-Sanfelice-Teel-book}.  The application of an hybrid invariance principle in the context of distributed event-based control requires some extra care, but it is rewarding and proves itself to be a powerful analytical tool. In this respect we view this as an additional contribution of the paper. We refer the reader to \cite{mayhew.et.al.tac2012} for other applications of hybrid stability tools for multi-agent cooperation.

Our results are applicable to systems subject to input saturation, which is also new when compared with existing event-based control results. We thus present simulation results for a network of two-dimensional linear systems subject to input saturations. Our preliminary work in \cite{De-Persis-Postoyan-mtns14} was dedicated to the rendez-vous for these particular systems in the case where the network is only composed of $2$ agents.  Compared to \cite{Liuzza-et-al-nolcos2013}, we address a different class of nonlinear systems as well as more general coordination tasks and we design time-triggered and self-triggered controllers based on a Lyapunov redesign.

The paper is organised as follows. Notations and preliminaries about the hybrid formalism of \cite{Goebel-Sanfelice-Teel-book} are provided in Section \ref{sect-prel}. The problem is stated in Section \ref{sect-pb-statement} and the event-triggered control strategies are presented in Section \ref{sect-etc}. The time-triggered and the self-triggered controllers are respectively developed in Sections \ref{sect-ttc} and \ref{sect-stc}. Section \ref{sect-simus} proposes simulations results. The proof of the main theorem is detailed in Section \ref{sect-proof-thm-etc}. Finally, Section \ref{sect-conclusions} concludes the paper.

\section{Preliminaries}\label{sect-prel}

Let $\Rl{}{}=(-\infty,\infty)$, $\Rlo=[0,\infty)$,
$\Rlp=(0,\infty)$, $\Zo=\{0,1,2,\ldots\}$, $\Zp=\{1,2,\ldots\}$. For $(x,y)\in\Rl{n+m}{}$,  $(x,y)$ stands for $[x^{\mathrm{T}},\,y^{\mathrm{T}}]^{\mathrm{T}}$. Let $f:\Rl{n}{}\rightarrow\Rl{}{}$ and $r\in\Rl{}{}$, we denote by $f^{-1}(r)$ the set $\{x\in\Rl{n}{}\,:\,f(x)=r\}$. 
A function $\gamma:\Rlo\rightarrow\Rlo$ is of class $\K$ if it is continuous, zero at zero and strictly increasing and it is of class $\Kinf$ if, in addition, it is unbounded. 
A set-valued mapping $M:\Rl{m}{}\rightrightarrows\Rl{n}{}$ is outer semicontinuous if and only if its graph $\{(x,y)\,:\, y\in M(x)\}$ is closed (see Lemma 5.10 in \cite{Goebel-Sanfelice-Teel-book}). The notation $\mathbb{I}$ denotes the identity matrix or application, and $\mathbf{1}$ and $\mathbf{0}$ are respectively the vector composed of $1$ and $0$ whose dimensions depend on the context. We use $\diag \{a_{1},\ldots,a_{n}\}$ to represent the diagonal matrix with constants $a_{1},\ldots,a_{n}$ on the diagonal. The Kronecker product of two matrices $A=[a_{ij}]\in\Rl{m\times n}{}$ and $B\in\Rl{p\times q}{}$ is written as
\[
A\otimes B= {\small\left[
                                                                                                                    \begin{array}{ccc}
                                                                                                                      a_{11}B & \ldots & a_{1n}B \\
                                                                                                                      \vdots & \ddots & \vdots \\
                                                                                                                      a_{m1}B & \ldots & a_{mn}B \\
                                                                                                                    \end{array}
                                                                                                                  \right]
}.
\]
We denote the distance of a point $x\in\Rl{n}{}$ to a set $\mathcal{A}\subset\Rl{n}{}$ as $\|x\|_{\mathcal{A}}=\inf\{\|x-y\|\,:\,y\in\mathcal{A}\}$. We recall the definition of the tangent cone to a set at a point (see Definition 5.12 in \cite{Goebel-Sanfelice-Teel-book}).

\begin{defn}\label{def-tangent-cone} The tangent cone to a set $S\subset\Rl{n}{}$ at a point $x\in\Rl{n}{}$, denoted
$T_{S}(x)$, is the set of all vectors $w\in\Rl{n}{}$ for which there exist $x_{i}\in S$, $\tau_{i}>0$ with
$x_{i}\to x$, $\tau_{i}\to 0$ as $i\to\infty$ such that $w=\lim_{i\to\infty}(x_{i}-x)/\tau_{i}$.  \hfill $\Box$
\end{defn}

We will study hybrid systems of the form below using the formalism of \cite{Goebel-Sanfelice-Teel-book}
\begin{equation}
\begin{array}{llllllllllllllll}
\dot{x} \in F(x) & \text{ for }x\in C, & &
x^{+} \in G(x) & \text{ for }x\in D,
\end{array}\label{eq-sys-prel}
\end{equation}
where $x\in\Rl{n}{}$ is the state, $F$ is the flow map, $G$ is the jump map, $C$ is the flow set and $D$ is the jump set.
We recall some definitions related to \cite{Goebel-Sanfelice-Teel-book}.
A subset $E\subset\Rlo\times\Zo$ is a \emph{hybrid time domain} if for all $(T,K)\in E$, $E\cap([0,T]\times\{0,\ldots,K\})={\underset{k\in\{0,1,\ldots,K-1\}}\bigcup}([t_{k},t_{k+1}],k)$ for
some finite sequence of times $0=t_{0}\leq t_{1} \leq \ldots \leq t_{K}$.
A function $\phi:E\rightarrow\Rl{n}{}$ is a \emph{hybrid arc} if $E$ is a hybrid time domain and if for each
$k\in\Zo$, $t\mapsto\phi(t,k)$ is locally absolutely continuous on $I^{k}=\{t\,:\,(t,k)\in E\}$.
We assume that: (i) $C$ and $D$ are closed subsets of $\Rl{n}{}$; (ii) $F$ is defined on $C$, is outer semicontinuous and locally bounded relative to $C$, and $F(x)$ is convex for every $x\in C$; (iii) $G$ is defined on $D$, is outer semicontinuous and locally bounded relative to $D$. The hybrid arc $\phi:\dom \phi\to \Rl{n}{}$ is a \emph{solution} to (\ref{eq-sys-prel}) if: (i) $\phi(0,0)\in C \cup D$;
(ii) for any $k\in\Zo$, $\phi(t,k)\in C$ and  $\frac{d}{dt}\phi(t,k)\in F(\phi(t,k))$ for almost all $t\in I^{k}$ (recall that $I^{k}=\{t\,:\,(t,k)\in \dom\phi\}$);
(iii) for every $(t,k)\in\dom \phi$ such that $(t,k+1)\in\dom \phi$, $\phi(t,k)\in D$ and $\phi(t,k+1)\in G(\phi(t,k))$.
A solution $\phi$ to (\ref{eq-sys-prel}) is: {\em nontrivial} if  $\dom \phi$ contains at least two points; \emph{maximal} if it cannot be extended; \emph{complete} if $\dom \phi$ is unbounded;
\emph{precompact} if it is complete and  the closure of its range is compact, where the range of $\phi$ is $\rge\phi:=\{y\in\Rl{n}{}\,:\, \exists (t,k)\in\dom \phi \text{ such that } y=\phi(t,k)\}$.

We introduce the following definition to denote solutions which have uniform average dwell-times.
\begin{defn}\label{def-prel-adt} The solutions to (\ref{eq-sys-prel}) have \emph{a uniform semiglobal average dwell-time} if for any $\Delta\geq 0$, there exist $\tau(\Delta)>0$ and $n_{0}(\Delta)\in\Zp$ such that for any solution $\phi$ to (\ref{eq-sys-prel}) with $\|\phi(0,0)\|\leq\Delta$
\begin{equation}
\begin{array}{rllll}
k-i & \leq & \dst\frac{1}{\tau(\Delta)}(t-s)+n_{0}(\Delta),
\end{array}\label{eq-def-adt}
\end{equation}
for any $(s,i),(t,k)\in\dom \phi$ with $s+i\leq t+k$. We say that the solutions to (\ref{eq-sys-prel}) have a \emph{uniform global average dwell-time} when $\tau$ and $n_{0}$ are independent of $\Delta$. \hfill $\Box$
\end{defn}

We recall the following invariance definition (see Definition 6.19 in \cite{Goebel-Sanfelice-Teel-book}).
\begin{defn}\label{def-fi}
A set $S\subset\Rl{n}{}$ is \emph{weakly invariant} for system (\ref{eq-sys-prel}) if it is:
\begin{itemize}
\item \emph{weakly forward invariant}, \ie for any $\xi\in S$ there exists at least one complete solution $\phi$ with initial condition $\xi$ such that $\rge \phi \subset S$;
\item \emph{weakly backward invariant}, \ie for any $\xi\in S$ and $\tau>0$, there exists at least one solution $\phi$ such that for some $(t^{*},k^{*})\in\dom\phi$, $t^{*}+k^{*}\geq \tau$, it is the case that $\phi(t^{*},k^{*})=\xi$ and $\phi(t,k)\in S$ for all $(t,k)\in\dom\phi$ with $t+k\leq t^{*}+k^{*}$. \hfill $\Box$
\end{itemize}
\end{defn}

Finally, we say that a solution $\phi$ \emph{approaches the set} $S\subset \Rl{n}{}$ (\cite{Sanfelice-Goebel-Teel-TAC07(invariance)}) if for any $\epsilon>0$ there exists $(t^{*},k^{*})\in\dom\phi$ such that for all $(t,k)\in\dom\phi$ with $t+k\geq t^{*}+k^{*}$, $\phi(t,k)\in S + \epsilon \mathbb{B}$, where $\mathbb{B}$ is the unit ball.


\section{Problem statement}\label{sect-pb-statement}

Our objective is to construct distributed controllers to ensure the coordination of networked systems with limited communication and/or computation resources. In particular, we consider $N$ agents which are interconnected over a connected\footnote{A graph is \emph{connected} if, for each pair of nodes $i,j$, there exists a path which connects $i$ and
$j$, where a path is an ordered list of edges such that the
head of each edge is equal to the tail of the following one.} undirected graph $\mathcal{G}=(\mathcal{I},\mathcal{E})$ where $\mathcal{I}:=\{1,\ldots,N\}$ is the set of nodes and $\mathcal{E}$ is the set of
pairs of nodes connected by edges.
The dynamics of the agents is given by
\begin{equation}
\begin{array}{lllll}
\dot{p}_{i} & = & y_{i}\\
\dot{v}_{i} & = & f_i(v_{i},u_{i})\\
y_i & = & h_i(v_i),
\end{array}\label{eq-sys-ct}
\end{equation}
where $p_{i}\in\Rl{n_p}{}$ and $v_{i}\in\Rl{n_{v_i}}{}$ are the states, $y_i\in\Rl{n_p}{}$ is the output, $u_{i}\in\Rl{n_p}{}$ is the control input, $f_i$ and $h_i$ are locally Lipschitz functions such that $h_i(0,0)=0$, and $f_i(0,u_i)=0$ implies that $u_i=0$, $i\in \mathcal{I}$. We note that the dimension of $v_i$ is agent-dependent and that the agents dynamics may be different,  hence the networked system is allowed to be heterogenous. Dynamical systems of the form of (\ref{eq-sys-ct}) can describe mechanical systems and vehicles (in which case $p_i$ and $v_i$ are typically the position and the velocity, respectively), as well as electrical devices to mention a few examples. To formally state our coordination goal, we need to introduce the \emph{relative distance}, for any $(i,j)\in\mathcal{E}$,
\begin{equation}
\begin{array}{lllll}
z_{ij} & := & p_{j}-p_{i}.
\end{array}\label{eq-z-ij}
\end{equation}
We want to ensure the convergence of every $z_{ij}$, $(i,j)\in\mathcal{E}$, to a prescribed compact set $\mathcal{A}_{ij}\subset\Rl{n_p}{}$, with $\mathcal{A}_{ij}=\mathcal{A}_{ji}$, as in \cite{Arcak-tac07}. The sets $\mathcal{A}_{ij}$ can be the origin, in which case the objective is to ensure the agreement among the agents' variables $p_i$'s, or it can be a vector different from the origin,  in which case we achieve a formation control, to give a few examples.

We follow an emulation approach to design the controllers. We first design the feedback laws $u_i$, $i\in\mathcal{I}$, in the ideal case where the agents have unlimited resources using the results of \cite{Arcak-tac07}. Afterwards,  we take into account the resources constraints to which are subject the agents and we synthesize appropriate triggering strategies to preserve the desired coordination task in this context. Since we design the feedback laws using \cite{Arcak-tac07}, we need to make the following assumption on the $v_i$-system, $i\in\mathcal{I}$.

\begin{ass}\label{ass-strict-passivity} For any $i\in\mathcal{I}$, the system $\dot v_{i} = f_{i}(v_i,u_i)$ is strictly passive from $u_i$ to $y_i=h_i(v_i)$ with a continuously differentiable storage function $S_i:\Rl{n_{v_i}}{}\to \Rlo$ such that there exist $\underline{\alpha}_{S_i},\overline{\alpha}_{S_i}\in\Kinf$, and  a positive definite function $\rho_i:\Rl{n_{v_i}}{}\to\Rlo$ which verify for any $v_i\in\Rl{n_{v_i}}{},u_i\in\Rl{n_{p}}{},y_i\in\Rl{n_p}{}$
\begin{equation}
\left\{\begin{array}{lllll}
\underline{\alpha}_{S_i}(\|v_i\|) \leq S_{i}(v_i) \leq \overline{\alpha}_{S_i}(\|v_i\|)\\
\left\langle \nabla S_i(v_i),f_i(v_i,u_i)\right\rangle \leq -\rho_i(v_i) + u_i^{\mathrm{T}}y_i.
\end{array}\right.\label{eq-ass-strict-passivity}
\end{equation} \hfill $\Box$
\end{ass}

Systems that satisfy Assumption \ref{ass-strict-passivity} have been widely investigated in the context of coordinating systems and appears in several applications (\cite{Arcak-tac07,Burger2013,schaft.maschke.sicon13}).
In continuous-time, the control input $u_{i}$ is defined as (\cite{Arcak-tac07})
\begin{equation}
\begin{array}{lllll}
u_{i} & = & \dst{\underset{j\in \Ni}\sum}\psi_{ij}(z_{ij})
\end{array}\label{eq-cont-ct}
\end{equation}
where $\Ni$ is the set of neighbours of the node $i\in \mathcal{I}$, \ie $\Ni:=\{j\in\mathcal{I}\,:\,(i,j)\in\mathcal{E}\}$.
The functions $\psi_{ij}:\Rl{n_p}{}\rightarrow\Rl{n_p}{}$, $(i,j)\in\mathcal{E}$, are designed as $\psi_{ij}=\nabla P_{ij}$ where $\nabla P_{ij}$ is the gradient of the designed function $P_{ij}:\Rl{n_p}{}\to\Rlo$ which is required to satisfy the following properties:
\begin{enumerate}
\item[(a)] $P_{ij}$ is is twice continuously differentiable;
\item[(b)] $P_{ij}=P_{ji}$;
\item[(c)] There exist $\underline{\alpha}_{P_{ij}},\overline{\alpha}_{P_{ij}}\in\Kinf$ such that
 $\underline{\alpha}_{P_{ij}}(\|x\|_{\mathcal{A}_{ij}}) \leq P_{ij}(x) \leq \overline{\alpha}_{P_{ij}}(\|x\|_{\mathcal{A}_{ij}})$ for any $x\in\Rl{n_i}{}$;
\item[(d)] $\psi_{ij}(-x)=-\psi_{ij}(x)$ for any $x\in\Rl{n_p}{}$.
\end{enumerate}
According to \cite{Arcak-tac07}, the controllers in (\ref{eq-cont-ct}) guarantee that, for any $(i,j)\in\mathcal{E}$, the relative distance $z_{ij}$ approaches the set $\mathcal{A}_{ij}$  (under an extra assumption specified later), which means that the coordination is achieved.

In this paper, we take into account the resources limitations of the system in terms of communication and/or computation. In particular, we envision a setting where the agents only receive measurements from their neighbours and/or update their control inputs at some given time instants to be determined. In this case, we denote the control input in (\ref{eq-cont-ct}) as $\hat u_i$ which is defined by, for $i\in \mathcal{I}$,
\begin{equation}
\begin{array}{lllll}
\hat u_{i} & = & \dst{\underset{j\in \Ni}\sum}\psi_{ij}(\hat{z}_{ij})
\end{array}\label{eq-cont}
\end{equation}
where $\hat{z}_{ij}$ is a sampled version of $z_{ij}$, which is locally maintained by agent $i$. 
This variable is held constant between two successive updates, \ie
$\dot{\hat{z}}_{ij} = 0$
and is reset to the actual value of $z_{ij}$ at the update time instant, which leads to the jump equation
\begin{equation}
\begin{array}{lllllllll}
\hat{z}_{ij}^{+} & = & z_{ij}.
\end{array}\label{eq-hat-z-ij-jump}
\end{equation}
A sequence of update time instants will be assigned to each pair $(i,j)\in \mathcal{E}$. These are time instants that are generated at agent $i$ and that are triggered by measurements relative to neighbor $j\in \mathcal{N}_i$. Symmetrically, agent $j$ will generate update time instants based on measurements relative to $i$.
The triggering conditions will be such that the  events generated by agent $i$ relative to neighbor $j$ and by agent $j$ relative to neighbor $i$ are the same. For this reason we term these instants as \emph{edge events}.
At each event of the edge $(i,j)\in \mathcal{E}$, the agents $i$ and $j$ communicate with each other and both of
them update the sampled variables $\hat{z}_{ij}$ and $\hat{z}_{ji}$ according to (\ref{eq-hat-z-ij-jump}), which leads to an update of the control inputs $\hat u_{i}$ and $\hat u_{j}$ in view of (\ref{eq-cont}).

Our goal is to define the sequence of edge events in order to save resources while still ensuring the desired coordination. We present solutions for the three scenarios listed below.
\begin{itemize}
\item \emph{Event-triggered control:} any agent knows its relative distance with any of its neighbours at any time instant and the corresponding part of the control input is only updated whenever a certain edge-dependent triggering condition is satisfied. 
    This setup requires that the agents are equipped with local sensors which measure the relative distance with their neighbour(s)  at a high frequency or that the agents communicate with their neighbour via a high-bandwidth communication channel. In that way, we can make the approximation that the agents continuously have access to their neighbour relative distance.
\item \emph{Time-triggered control:} any agent has access to its relative distances and updates its control input \emph{only} at edge-dependent time instants which are generated by a time-driven policy. These edge events can be periodic, but that is not necessary: we do allow aperiodic sampling.
\item \emph{Self-triggered control:} any agent has access to the relative distance as well as its time derivative and updates the corresponding sampled variables \emph{only} at edge events. 
    The next edge event is determined by the values of the relative distance and its time derivative at the last transmission. This scheme reduces the usage of the agents sensors or of the communication channel, and potentially of the agent CPU, as we will explain later. It typically generates more edge events compared to event-triggered control (but it does not require the continuous measurement of the neighbours relative distance) and less events than time-triggered control, see for example the simulation results in Section \ref{sect-simus}.
\end{itemize}
The proposed strategies ensure the existence of a uniform strictly positive amount of time between two successive events of a given edge. 
We do tolerate the occurrence of a finite number of simultaneous edge events for a given agent as in \emph{e.g.}, \cite{Dimarogonas-et-al-tac2012,De-Persis-Frasca-TAC'13}. We assume that the agent hardware handles this situation by prioritizing the edge events, which typically leads to small-delays in the control input. We do not address the analysis of the effect of these delays in this paper.

\begin{rem}\label{rem-velocity-requirement} We have not specified any requirement on the states $v_i$, $i\in\mathcal{I}$, for the coordination objective. We will see in the next sections that these variables converge to the origin. The extension to the case where $v_i$ has to converge to a prescribed time-varying vector $\mathbf{v}_i$ as in \cite{Arcak-tac07} is left for future work. The reason is the following. In a realistic setting, only a sampled version of $\mathbf{v}_i$ can be available to the agent $i\in\mathcal{I}$. This sampling typically generates errors which affect the asymptotic convergence of $v_i$ to $\mathbf{v}_i$ and leads to technical difficulties, as shown in \cite{Postoyan-van-de-Wouw-et-al-paper} in the context of networked control systems. Note though that our results directly apply when the $\mathbf{v}_i$'s are constant. In this case, following \cite{Arcak-tac07},  $\dot p_i = y_i + \mathbf{v}_i$ in (\ref{eq-sys-ct}), and only one sample is needed to generate $\mathbf{v}_i$ since the latter takes a constant value.
\hfill $\Box$
\end{rem}

\section{Event-triggered control}\label{sect-etc}

\subsection{Triggering conditions and hybrid model}\label{subsect-hybrid-model}

Consider the agent $i\in \mathcal{I}$. To define the events associated with the edge $(i,j)$ where $j\in\Ni$, we introduce an auxiliary variable $\phi_{ij}\in\Rl{}{}$, which we call a \emph{clock}. The idea is to reset $\phi_{ij}$ to a constant value $b_{ij}>0$ after each event associated with $(i,j)$ and to trigger the next one when $\phi_{ij}$ becomes equal to $a_{ij}\in[0,b_{ij})$.  The constants $a_{ij}$ and $b_{ij}$ are designed parameters. Between two successive edge events, $\phi_{ij}$ is given by the solution to the ordinary differential equation below
\begin{equation}
\begin{array}{lllll}
\dot{\phi}_{ij} & = & \dst -\frac{1}{\sigma_{ij}}\left(1+\phi_{ij}^{2}\left\|\nabla \psi_{ij}(z_{ij})\right\|^{2}\right),
\end{array}\label{eq-phi-ij}
\end{equation}
where
$\sigma_{ij}$ is a strictly positive constant which will be specified in the following,
$\|\nabla \psi_{ij}(z_{ij})\|$ is the induced matrix Euclidean norm of the matrix $\nabla \psi_{ij}(z_{ij})$,
and we recall that $z_{ij}=p_{j}-p_{i}$. We notice that $\phi_{ij}$ strictly decreases on flows in view of (\ref{eq-phi-ij}). The length of the inter-event times depends on the choice of the constants $a_{ij}$ and $b_{ij}$. To take $a_{ij}$ small and $b_{ij}$ large typically helps enlarging the inter-event time, at the price of a degraded speed of convergence as the evolution of the variables $v_i$ depends on the sampled control input, see for an illustration the simulation results in Section \ref{sect-simus}.
The clock $\phi_{ij}$ can be locally implemented on agent $i$ provided that continuous measurements of $z_{ij}$ are available, which is assumed to be the case in this section.

\begin{remark}
The clock dynamics (\ref{eq-phi-ij}) descends from the Lyapunov analysis carried out in Section  \ref{subsect-proof-lyap}. To help the reader grasping the significance of (\ref{eq-phi-ij}), we provide here a preliminary discussion. In Section  \ref{subsect-proof-lyap}, we first introduce an energy-like Lyapunov function which is commonly used  in the stability analysis of the networked systems (\ref{eq-sys-ct}), see \cite{Arcak-tac07}. Then we show that during the continuous evolution of (\ref{eq-sys-ct}) under the sampled-data control (\ref{eq-cont}) (see (\ref{eq-sys-etc}) below for a formal description of the overall dynamical system under consideration), if the sampling occurs according to rule (\ref{eq-phi-ij}), then the energy-like function extended to include the `energy' associated with the sampling errors is monotonically non-increasing. \hfill $\Box$
\end{remark}

The dynamics of the agent $i\in \mathcal{I}$ can be described by the hybrid system below
\begin{equation}
\begin{array}{l}
\left.\begin{array}{llllll}
\dot{p}_{i} & = & y_{i}\\
\dot{v}_{i} & = & f_i(v_{i}, \hat u_i)\\
\dot{\hat{z}}_{ij} &=& 0& j\in \Ni\\
\dot{\phi}_{ij} & = & -\frac{1}{\sigma_{ij}}\left(1+\phi_{ij}^{2}\left\|\nabla \psi_{ij}(z_{ij})\right\|^{2}\right)
& j\in\Ni
\end{array}\right\}  \\[4mm ] \hfill  \forall j\in \Ni\,\,\,\phi_{ij}\in[a_{ij},b_{ij}]\\[4mm ]
\left.
\begin{array}{rlllll}
p_{i}^{+} & = & p_{i}\\
v_{i}^{+} & = & v_{i}\\
\left(
  \begin{array}{c}
    \hat{z}_{ij}^+ \\
    \phi_{ij}^{+} \\
  \end{array}
\right)  & = & \!\!\!\!\left\{\begin{array}{lllll}
                 \begin{array}{llllllll}
                   \left(
                      \begin{array}{c}
                        z_{ij} \\
                        b_{ij} \\
                      \end{array}\right) & j\in\Ni \text{ and } \phi_{ij}=a_{ij}  \\
                 \left(
                \begin{array}{llllllllllllllll}
                    \hat{z}_{ij} \\
                    \phi_{ij} \\
                \end{array}\right) & j\in\Ni \text{ and } \phi_{ij} > a_{ij} \\
                 \end{array}
               \end{array}\right.
\end{array}\!\!\!\!\right\}  \\[4mm ]
\hfill   \exists j\in \Ni\,\,\,\phi_{ij}=a_{ij},
\end{array}
\label{eq-sys-etc-i}
\end{equation}
where $\hat u_i$ is defined in (\ref{eq-cont}). The jump map in (\ref{eq-sys-etc-i}) means that only the pairs $(\hat{z}_{ij},\phi_{ij})$, $j\in\Ni$, for which $\phi_{ij}$ is equal to $a_{ij}$, are reset to $(z_{ij},b_{ij})$; the others remain unchanged. We see that the control input updates are edge-dependent and distributed as desired.
In the analysis that follows, it is essential that each agent $i$ maintains a local sampled version of the measurement $z_{ij}$, $j\in \Ni$, which is consistent with the local sampled version of the corresponding quantity $z_{ji}$ by the agent $j$. To be more specific, for $(i,j)\in \mathcal{E}$, it must be true that  $\hat{z}_{ij}(t,k)=-\hat{z}_{ji}(t,k)$ for all $(t,k)$ in the domain of the solution. To guarantee this property, we make the following assumption.

\begin{ass}\label{ass-synchro} The following hold for any $(i,j)\in \mathcal{E}$.
\begin{itemize}
\item[(i)] $a_{ij}=a_{ji}$, $b_{ij}=b_{ji}$, $\sigma_{ij}=\sigma_{ji}$.
\item[(ii)] The variables $\hat{z}_{ij}$ and $\phi_{ij}$ are respectively initialized at the same values as $-\hat{z}_{ji}$ and $\phi_{ji}$.
\hfill $\Box$
\end{itemize}
\end{ass}

Assumption \ref{ass-synchro} introduces no major conservatism as neighboring agents can {\em a priori} agree on the constants $a_{ij},\,a_{ji},\,b_{ij},\,b_{ji},\,\sigma_{ij},\,\sigma_{ji}$ and the initial conditions $\hat{z}_{ij}$ and $\phi_{ij}$. Notice in particular that, in the analysis below, the  initial condition for $\hat{z}_{ij}$ must not necessarily be set equal to the measured quantity ${z}_{ij}$.
When Assumption \ref{ass-synchro} is not verified, the clocks $\phi_{ij}$ and $\phi_{ji}$, $(i,j)\in \mathcal{E}$, will be
different and this will imply that the updates for $\hat z_{ij}$ and $\hat z_{ji}$ will occur at different times and that the two measurements are different. This causes an asymmetry in the control laws of the neighboring agents $i, j$ that may disrupt the convergence of the algorithms. Robustness of our algorithm to asymmetric initializations is an important open problem.

\begin{rem}
In different scenarios, item (ii) of Assumption \ref{ass-synchro} may be less critical. In fact, the scenario that was discussed above assumes that when the clock $\phi_{ij}$ reaches $a_{ij}$, the agent  updates $\hat{z}_{ij}$ with the information collected by its sensor. A different scenario could be as follows. Assume that the two clock variables $\phi_{ij}$ and $\phi_{ji}$, $(i,j)\in \mathcal{E}$, are initially different until one of these, say $\phi_{ij}$, becomes equal to $b_{ij}$ (recall that $b_{ij}=b_{ji}$ in view of item (i) of Assumption \ref{ass-synchro}).
At this time instant, we can envision the case in which agent $i$ (the one whose clock variable has become equal to $b_{ij}$) notifies (without delay) agent $j$ to update its own clock variable.
Hence, $(\hat{z}_{ij},\phi_{ij})$ and $(\hat{z}_{ji},\phi_{ji})$ are updated respectively to $(z_{ij},b_{ij})$ and $(z_{ji},b_{ij})$. In that way, the pairs $(\phi_{ij},\hat{z}_{ij})$ and $(\phi_{ji},-\hat{z}_{ji})$ are equal for all future times in view of\footnote{Note that $\left(\nabla \psi_{ij}(z_{ij})\right)^{2}=\left(\nabla \psi_{ji}(z_{ji})\right)^{2}$ in (\ref{eq-phi-ij}) as $z_{ij}=-z_{ji}$ from (\ref{eq-z-ij}) and since $\psi_{ij}$ satisfies item (d) in  Section \ref{sect-pb-statement}.} (\ref{eq-sys-etc-i}) and the convergence results presented hereafter do apply in this case. \hfill $\Box$
\end{rem}

In view of Assumption \ref{ass-synchro}, we no longer need to distinguish $\phi_{ij}$ from $\phi_{ji}$. We can therefore define a single clock $\phi_{\ell}$ instead, where $\ell$ is the index associated with the edge $(i,j)\in\mathcal{E}$. A similar remark applies for the sampled variables $\hat{z}_{ij}$ and $\hat{z}_{ji}$ as $\hat{z}_{ij}=-\hat{z}_{ji}$. For that purpose, we assign to each edge of $\mathcal{E}$ an arbitrary direction and we denote by $M$ the number of edges of the graph $\mathcal{G}$ which we number. We define the incidence matrix $D$ of $\mathcal{G}$ as $D=[d_{i\ell}]_{(i,\ell)\in \mathcal{I}\times\{1,\ldots,M\}}$ with
$d_{i\ell}=1$ if the node $i$ is the positive end of the $\ell^{\text{th}}$ edge, $d_{i\ell}=-1$ if the agent $i$ is the negative end of the $\ell^{\text{th}}$ edge, and $d_{i\ell}=0$ otherwise.
In that way, we define, for the $\ell^{\text{th}}$ edge corresponding to $(i,j)\in\mathcal{E}$,
\[
z_\ell:=\left\{
\begin{array}{cl}
z_{ij} & \textrm{if $j$ is the positive end of the edge $\ell$}\\
z_{ji} & \textrm{if $i$ is the positive end of the edge $\ell$},
\end{array}
\right.
\]
and
\[
\hat{z}_\ell:=\left\{
\begin{array}{cl}
\hat{z}_{ij} & \textrm{if $j$ is the positive end of the edge $\ell$}\\
\hat{z}_{ji} & \textrm{if $i$ is the positive end of the edge $\ell$}.
\end{array}
\right.
\]
For the $\ell^{\text{th}}$ edge corresponding to $(i,j)\in\mathcal{E}$, we rewrite the dynamics in (\ref{eq-phi-ij}) as
\begin{equation}
\begin{array}{lllll}
\dot{\phi}_{\ell} & = & \dst-\frac{1}{\sigma_{\ell}}\left(1+\phi_{\ell}^{2}\left\|\nabla \psi_{\ell}(z_{\ell})\right\|^{2}\right),
\end{array}\label{eq-phi-l}
\end{equation}
where  $\sigma_{\ell}:=\sigma_{ij}=\sigma_{ji}$, $a_{\ell}:=a_{ij}=a_{ji}$ and $b_{\ell}:=b_{ij}=b_{ji}$ (in view of Assumption \ref{ass-synchro}). We similarly define $\mathcal{A}_{\ell}:=\mathcal{A}_{ij}=\mathcal{A}_{ji}$ and $P_{\ell}=P_{ij}$ where $(i,j)\in\mathcal{E}$ is the $\ell^{\text{th}}$ edge.

We are not ready yet to present a model of the overall system. Indeed, it appears that the map which defines the jump equation in (\ref{eq-sys-etc-i}) and which becomes with the notation introduced above, with $\mathcal{E}_{i}$ the set of edge indices corresponding to the  edges connected to agent $i$,
\begin{equation}
\begin{array}{rlllll}
p_{i}^{+} & = & p_{i}\\
v_{i}^{+} & = & v_{i}\\
\left(
  \begin{array}{c}
    \hat{z}_{\ell}^+ \\
    \phi_{\ell}^{+} \\
  \end{array}
\right)  & = & \left\{\begin{array}{lllll}
                 \begin{array}{llllllll}
                   \left(
                      \begin{array}{c}
                        z_{\ell} \\
                        b_{\ell} \\
                      \end{array}\right) & \ell\in  \mathcal{E}_i \text{ and } \phi_{\ell}=a_{\ell}  \\
                 \left(
                \begin{array}{llllllllllllllll}
                    \hat{z}_{\ell} \\
                    \phi_{\ell} \\
                \end{array}\right) & \ell\in \mathcal{E}_i  \text{ and } \phi_{\ell} > a_{\ell}, \\
                 \end{array}
               \end{array}\right.
\end{array}\label{eq-sys-etc-l}
\end{equation}
is not outer semicontinuous because its graph is not closed. This is an issue because the outer semicontinuity of the jump map is a necessary condition for a hybrid system to be (nominally) well-posed (see Lemma 6.9 in \cite{Goebel-Sanfelice-Teel-book}) which is required to apply the invariance principles presented in Chapter 8 in \cite{Goebel-Sanfelice-Teel-book}.

To overcome that issue, we redefine the jump map. We use the technique proposed in \cite{Sanfelice-Biemond-et-al-ijrnc13} for that purpose. Instead of doing it for the model of a single agent, we directly do it on a model of the overall system. Hence, we define the concatenated vectors $p:=(p_{1},p_{2},\ldots,p_{N})\in\Rl{n_p N}{}$, $v:=(v_{1},v_{2},\ldots,v_{N})\in\Rl{n_v}{}$, $\hat u:=(\hat u_{1},\hat u_{2},\ldots,\hat u_{N})\in\Rl{n_p N}{}$, $\phi:=(\phi_{1},\ldots,\phi_{M})\in\Rl{M}{}$, $z:= (z_{1},\ldots,z_{M})\in\Rl{n_p M}{}$, and $\hat{z}:=(\hat{z}_{1},\ldots,\hat{z}_{M})\in\Rl{n_p M}{}$, with $n_v:=\sum_{i\in\mathcal{I}}n_{v_i}$. The system is modeled as follows
\begin{equation}
\begin{array}{rllll}
\left.\begin{array}{lllll}
\dot{p} & = & h(v) \\
\dot{v} & = & f(v,\hat u) \\
\dot{\hat{z}} &=& \mathbf{0}\\
\dot{\phi} & = & -\Sigma^{-1}
\left(\mathbf{1}+ \Phi(z)\right)\\
\end{array}\right\} &
\hspace{-0.2cm}
     \begin{array}{l}\forall \ell\in\{1,\ldots,M\}\\
                     \phi_{\ell}\in [a_{\ell},b_{\ell}]
     \end{array}\\
\left.
\begin{array}{rllll}
p^{+} & = & p\\
v^{+} & = & v\\
\left(\begin{array}{l}
\hat{z}^+\\
\phi^{+}
\end{array}\right) & \in &
G(z,\hat{z},\phi)
\end{array}\right\} & \hspace{-0.2cm}
     \begin{array}{l}\exists\ell\in\{1,\ldots,M\} \\
                     \phi_{\ell}=a_{\ell},
     \end{array}
\end{array}\label{eq-sys-etc}
\end{equation}
where $h(v):=(h_1(v_1),\ldots,h_N(v_N))$, $f(v,\hat u):=(f_1(v_1,\hat u_1),\ldots,f_N(v_N,\hat u_N))$, $\Sigma:=\text{diag}\{\sigma_1,\ldots,$ $\sigma_M\}$ and $\Phi(z):=(\phi_{1}^{2}\left\|\nabla \psi_1(z_1)\right\|^2,$ $\ldots,\phi_{M}^{2}\left\|\nabla \psi_M(z_M)\right\|^2)$.
Inspired by \cite{Sanfelice-Biemond-et-al-ijrnc13}, the set-valued jump map $G$ is defined as, for $(z,\hat{z},\phi)\in\Rl{(2n_p+1)M}{}$,
\begin{equation}
G(z,\hat{z},\phi)  :=  \left\{G_{\ell}(z,\hat{z},\phi)\,:\,\ell\in\{1,\ldots,M\}\,\text{ and }\,\phi_{\ell}=a_{\ell}\right\},
\label{eq-jump-map-revised}
\end{equation}
with, for $\ell\in\{1,\ldots,M\}$,
\begin{equation}
\begin{array}{lll}
G_{\ell}(z,\hat{z},\phi) &:=& \big(\hat{z}_{1},\ldots,\hat{z}_{\ell-1},z_{\ell},\hat{z}_{\ell+1},\ldots, \hat{z}_{M},\\
                          & & \quad\quad\quad  \phi_{1},\ldots, \phi_{\ell-1},b_{\ell},\phi_{\ell+1},\ldots,\phi_{M}\big).
\end{array}\label{eq-G-ell}
\end{equation}
In that way, when the clock $\phi_{\ell}$ is the only one which is equal to its lower bound $a_{\ell}$, the pair $(\phi_{\ell},\hat{z}_{\ell})$ is reset to $(b_{\ell},z_{\ell})$, while the others remain unchanged.
In contrast to (\ref{eq-sys-etc-l}), when several clocks have reached their lower bound, the jump map (\ref{eq-jump-map-revised}) only allows a single edge to reset its clock and its sampled variable. Consequently, a finite number of jumps successively occurs in this case (with no flow in between), until all the concerned edge variables have been updated.
A couple of remarks about system (\ref{eq-sys-etc}) need to be added. First, the map $G$ in (\ref{eq-jump-map-revised}) is defined on $\Rl{(2n_p+1)M}{}$. When the states are in the jump set its definition is clear from (\ref{eq-jump-map-revised}), when these are not in the jump set, \ie when $\phi_{\ell}\neq a_{\ell}$ for any $\ell\in\{1,\ldots,M\}$, it reduces to the empty set. Second, $G$ is indeed outer semicontinuous as its graph is given by the union of the graphs of the mappings $G_{\ell}$, $\ell\in\{1,\ldots,M\}$, which are closed since these mappings are continuous. We also note that $G$ is locally bounded. As a consequence, since the flow map is continuous and the flow and the jump sets are closed, system (\ref{eq-sys-etc}) is (nominally) well-posed (see Theorem 6.30 in \cite{Goebel-Sanfelice-Teel-book}) and we will be able to apply the hybrid invariance principle in Chapter 8 of \cite{Goebel-Sanfelice-Teel-book} to investigate convergence.

\subsection{Main result}\label{subsect-main-result}

We are ready to state the main result of this section. The proof is provided in Section \ref{sect-proof-thm-etc}.
\begin{thm}\label{th-etc-N} Consider system (\ref{eq-sys-etc}) and suppose the following holds.
\begin{enumerate}
\item[(i)] Assumptions \ref{ass-strict-passivity}-\ref{ass-synchro} hold.
\item[(ii)] There exist $\kappa_{1},\ldots,\kappa_{M}\in(0,1)$  such that, for any $i\in\mathcal{I}$ and $v_i\in\Rl{n_{v_i}}{}$,
\begin{equation}
\begin{array}{lllll}
\dst -\kappa_i \rho_i(v_i) +2 \deg_i \max_{\ell \in \mathcal{E}_i}\{\sigma_\ell\}\left\|h_i(v_i)\right\|^2  & \le & 0
\end{array}\label{eq-thm-condition}
\end{equation}
where the $\sigma_\ell$'s come from (\ref{eq-phi-l}) and $\deg_i$ is the degree of agent $i$, \ie  the number of edges incident to agent $i$.
\item[(iii)] For any $z\in\Rl{n_p M}{}$, $(D \otimes \mathbb{I})\Psi(z)=\mathbf{0}$ implies $z\in\mathcal{A}$, where $\Psi(z):=(\psi_1(z_1),\ldots,\psi_M(z_M))$ and $\mathcal{A}:=\mathcal{A}_1\times\ldots\times\mathcal{A}_M$.
\end{enumerate}
The solutions have a uniform semiglobal average dwell-time and the maximal solutions are complete and approach the set $\left\{(p,v,\hat{z},\phi): z\in\mathcal{A},\,\,v=\mathbf{0},\,\phi_{\ell}\in[a_{\ell},b_{\ell}] \text{ for }\ell\in \{1,\ldots,M\}\right\}$. \hfill $\Box$
\end{thm}

Item (iii) in Theorem \ref{th-etc-N} is Assumption 1 in \cite{Arcak-tac07} (note that in our case $z$ always lies in the range space of $D^{\mathrm{T}}\otimes \mathbb{I}$ since the $P_\ell$'s are defined on $\Rl{n_p}{}$). In the proof of Theorem \ref{th-etc-N}, we show that $(D\otimes \mathbb{I})\Psi(z)$ converges to the origin, thus showing convergence of $z$ to the desired target set $\mathcal{A}$ in view of condition (iii).  The validity of this condition depends on the set $\mathcal{A}$. It is satisfied by important coordination tasks, such as rendez-vous and formation control (cf., \eg \cite{Arcak-tac07,Bai2011}).

We see that we need an extra condition to hold  compared to \cite{Arcak-tac07}, namely (\ref{eq-thm-condition}). It is satisfied when
\begin{equation}
\begin{array}{llllll}
\left\|h_{i}(v_i)\right\|^{2} & \leq &  C_i \rho_i(v_i) & & \forall v_i\in\Rl{n_{v_i}}{},
\end{array}\label{eq-thm-ass-cond}
\end{equation}
for some $C_i\in\Rlo$ and $i\in\mathcal{I}$. Indeed, it suffices to take $\sigma_\ell$, $\ell\in\mathcal{E}_i$ and $i\in\mathcal{I}$, sufficiently small such that, for a given $\kappa_i\in(0,1)$,
\begin{equation}
\begin{array}{llllll}
2 \dst\deg_i \max_{\ell \in \mathcal{E}_i}\{\sigma_\ell\} C_i& \leq  & \dst \kappa_i.
\end{array}\label{eq-thm-ass-param}
\end{equation}
Inequality (\ref{eq-thm-ass-param}) is equivalent to $C_i  \leq  \frac{\kappa_i}{2 \deg_i \max_{\ell \in \mathcal{E}_i}\{\sigma_\ell\}}$, which leads to $\left\|h_{i}(v_i)\right\|^{2} \leq  \frac{\kappa_i}{2 \deg_i \max_{\ell \in \mathcal{E}_i}\{\sigma_\ell\}}\rho_i(v_i)$ for any $v_i\in\Rl{n_{v_i}}{}$ in view of (\ref{eq-thm-ass-cond}), which in turn ensures (\ref{eq-thm-condition}). We notice that each agent only needs to know the degree of its neighbours and the local constant $C_i$ to synthesize its constants $\sigma_{\ell}$ in this case, $\ell\in\mathcal{E}_i$.  The knowledge of the agent degree can be achieved via an initial communication round during which the agents communicate their degrees to their neighbours.

\begin{rem} The fact that an additional condition is needed to prove the desired asymptotic convergence property under the considered sampling effects is in agreement with the literature on the stabilization of nonlinear sampled-data systems. Indeed, we know from \cite{Nesic_Teel_TAC04} that only semiglobal and practical stability can be ensured in general when emulating a globally asymptotically stabilizing continuous-time controller with fast sampling (under mild conditions); additional properties are needed to preserve asymptotic stability, like in Theorem \ref{th-etc-N}. \hfill $\Box$
\end{rem}

As mentioned in Section \ref{sect-pb-statement}, we cannot guarantee the existence of a dwell-time for the overall system as several agents may update their control inputs at the same instant or the same agents may have several of its local triggering conditions simultaneously violated. However, we do guarantee the existence of a uniform (semiglobal) dwell-time for each edge event (see Section \ref{subsect-proof-avg-dwt}), which in turn ensure the existence of a uniform semiglobal average dwell-time for the solutions of the overall system as stated in Theorem \ref{th-etc-N}.

\section{Time-triggered control}\label{sect-ttc}

In this section, we aim at defining the edge events using time-triggered rules. We rely for that purpose on the event-triggering strategies developed in the previous section which ensure the existence of a semiglobal dwell-time for each edge. In other words, there exists a strictly positive bound on the minimum time between two successive edge events, which depends on the ball of initial conditions (see Section \ref{subsect-proof-avg-dwt} for more details). We could use these dwell-times as an upper-bound on the \emph{maximum allowable time between two edge events} (MATE) to derive time-triggered strategies. However the fact that these constants depend on the ball of initial conditions render their implementation hard to achieve in practice, as each agent would need to know
the initial conditions of the other agents (more precisely the constant $\bar \Delta$ in Section \ref{subsect-proof-avg-dwt} which does depend on the agents' initial conditions) to compute its MATEs.
To overcome this issue, we design the function $\psi_{\ell}$, $\ell\in\{1,\ldots,M\}$, such that the following property holds, in addition to those listed in Section \ref{sect-pb-statement},
\begin{equation}
\forall \ell\in\{1,\ldots,M\}\,\,\,\exists K_\ell\geq 0\,\,\,\forall z_{\ell}\in\Rl{n_p}{} \,\,\, \|\nabla \psi_{\ell}(z_\ell) \| \leq K_\ell.
\label{eq-boundedness-nabla-psi}
\end{equation}
Property (\ref{eq-boundedness-nabla-psi}) is verified when $\psi_{\ell}$, $\ell\in\{1,\ldots,M\}$,  is globally Lipschitz. We denote the MATE of edge $\ell\in\{1,\ldots,M\}$ as $T_{\ell}$. The constant $T_{\ell}$ is the time it takes for the solution $\theta_\ell$ to the differential equation
\begin{equation}
\begin{array}{rllllllll}
\dot \theta_\ell = -\dst\frac{1}{\sigma_{\ell}} (1+\theta_\ell^2 K_{\ell}^2 ),  & &  \theta_\ell(0) = b_{\ell},
\end{array}\label{eq-ttc-theta}
\end{equation}
to decrease to $a_{\ell}$, like in \cite{Nesic_Teel_Carnevale_TAC09}. Equation (\ref{eq-ttc-theta}) corresponds to (\ref{eq-phi-l}) where $\|\nabla\psi_\ell(z_\ell)\|$ is replaced by its upper-bound $K_\ell$. In that way, the dynamics of $\theta_\ell$ is independent of the state. The solution to the differential equation given the initial condition $\theta_\ell(0)=b_{\ell}$ verifies, for $t\geq 0$,
$\arctan(K_\ell\theta_\ell(t)) = \arctan(K_\ell b_\ell)-\displaystyle\frac{K_\ell}{\sigma_\ell}t$,
from which it is inferred that
\begin{equation}
\begin{array}{rllllllll}
T_\ell & := & \displaystyle\frac{\sigma_\ell}{K_\ell}(\arctan(K_\ell b_\ell)-\arctan(K_\ell a_\ell)).
\end{array}\label{eq-T}
\end{equation}
Since $a_\ell,b_\ell$ can be chosen arbitrarily, the sampling interval can be changed, although it can never be larger than $\frac{\sigma_\ell}{K_\ell}\frac{\pi}{2}$ in view of (\ref{eq-T}). However, this choice might affect the speed of convergence of the system as the evolution of the velocities depends on the sampled control input. We represent the system using the hybrid model below, like in \cite{Nesic_Teel_Carnevale_TAC09},
\begin{equation}
\begin{array}{rllll}
\left.\begin{array}{lllll}
\dot{p} & = & h(v)\\
\dot{v} & = & f(v,\hat u) \\
\dot{\hat{z}} &=& \mathbf{0}\\
\dot{\tau} & = & \mathbf{1}
\end{array}\right\} &
\hspace{-0.2cm} \forall \ell\in\{1,\ldots,M\}\,\,\,\tau_{\ell}\in [0,T_{\ell}]\\
\left.
\begin{array}{rllll}
p^{+} & = & p\\
v^{+} & = & v\\
\left(\begin{array}{l}
\hat{z}^+\\
\tau^{+}
\end{array}\right) & \in &
\Gamma(z,\hat{z},\tau)
\end{array}\right\} & \hspace{-0.2cm}\begin{array}{l}\exists \ell\in\{1,\ldots,M\}\\
                                        \,\,\,\tau_{\ell}\in[\epsilon_{\ell},T_{\ell}],
                                     \end{array}
\end{array}\label{eq-sys-ttc}
\end{equation}
where $\tau:=(\tau_1,\ldots,\tau_M)\in\Rl{M}{}$ and $\tau_{\ell}$ is the time elapsed since the last event for the edge $\ell\in\{1,\ldots,M\}$. The constants $\epsilon_\ell$ can take any value in $(0,T_\ell]$ and represent the required minimum time between two successive events of edge $\ell$ to prevent arbitrarily close-in-time updates. This definition of the jump set allows to model the scenario where the edge events are not necessarily periodic but occur at most every $\epsilon_\ell$ units of times and at least every $T_{\ell}$ units of time. The function $\Gamma$ is defined in a similar way as $G$ in (\ref{eq-sys-etc})
\begin{equation}
\begin{array}{lll}
\Gamma(z,\hat{z},\tau) & := & \{\Gamma_{\ell}(z,\hat{z},\tau)\,:\\
                       &  & \quad \ell\in\{1,\ldots,M\}\,\text{ and }\,\tau_{\ell}\in[\epsilon_\ell,T_{\ell}]\},
\end{array}\label{eq-jump-map-revised-Gamma}
\end{equation}
with, for $\ell\in\{1,\ldots,M\}$,
$\Gamma_{\ell}(z,\hat{z},\tau) \!:= \!\big(\!\hat{z}_{1},\ldots,\hat{z}_{\ell-1},z_{\ell},$ $\hat{z}_{\ell+1},\ldots,\hat{z}_{M},
\!\tau_{1},\ldots, \!\tau_{\tau-1},0,\!\tau_{\ell+1},\ldots,\tau_{M}\!\big)$.

The result below follows from the proof of Theorem \ref{th-etc-N}.
\begin{cor}\label{cor-ttc-N} Consider system (\ref{eq-sys-ttc}) and suppose the following holds.
\begin{enumerate}
\item[(i)] Items (i)-(iii) of Theorem \ref{th-etc-N} hold.
\item[(ii)] Property (\ref{eq-boundedness-nabla-psi}) is guaranteed.
\end{enumerate}
The solutions have a uniform global average dwell-time and the maximal solutions are complete and approach the set
$\left\{(p,v,\hat{z},\tau): z\in\mathcal{A},\,v=\mathbf{0},\,\tau_{\ell}\in[0,T_{\ell}] \text{ for }\,\ell\in \{1,\ldots,M\}\right\}$. \hfill $\Box$ 
\end{cor}

Corollary \ref{cor-ttc-N} means that the variable $z$ is guaranteed to approach the prescribed compact set $\mathcal{A}$ as desired and the variable $v$ converges to the origin. The main difference with Theorem \ref{th-etc-N} is that a uniform global average dwell-time is guaranteed to exist, as opposed to a uniform semiglobal average dwell-time in Theorem \ref{th-etc-N}. This is possible due to the satisfaction of (\ref{eq-boundedness-nabla-psi}).

\section{Self-triggered control}\label{sect-stc}

The time-triggered implementation in the previous section is easy to implement but it has the drawback that the sampling at each edge $\ell\in\{1,\ldots,M\}$ is independent of the current value of $z_\ell$ and as such it might lead to some conservatism. On the other hand, the event-based control strategy of Section \ref{sect-etc} takes full advantage of $z_\ell$,  measuring it continuously over the inter-sampling period. Self-triggered control offers a compromise between these two paradigms. The idea is to define the MATE based on the values of the relative distance and its time derivative at the last edge event. In that way, the MATE is adapted to the current state of the system, as opposed to the time-triggered implementation, and the relative distance is not continuously monitored as in event-triggered control.
Recall that in event-triggered control, for each $\ell\in \{1,\ldots,M\}$, the sampling is dictated by the clock variable $\phi_\ell$  that flows according to $
\dot{\phi}_\ell = -\frac{1}{\sigma_{\ell}}\left(1+\phi_\ell^{2}\left\|\nabla \psi_{\ell}(z_\ell)\right\|^{2}\right)$.
To prevent the continuous measurement of $z_\ell$, the idea here is to replace $\left\|\nabla \psi_{\ell}(z_\ell)\right\|^{2}$ with a suitable function $\lambda_\ell$, which only depends on the value of $z_\ell$ and its time derivative at the last edge event.

\subsection{Construction of $\lambda_\ell$}
To preserve the properties ensured by the event-triggered controllers in Section \ref{sect-etc}, the function $\lambda_\ell$ has to be an upper bound on $\left\|\nabla \psi_{\ell}(z_\ell)\right\|^{2}$ (just like $K_\ell$ upper-bounds $\left\|\nabla \psi_{\ell}(z_\ell)\right\|$ in Section \ref{sect-ttc}). In that way, 
we will be able to apply the same arguments as for event-triggered control to analyse convergence. To derive such a bound, an estimate of $z_\ell$ is needed. As a matter of fact, if two vector-valued maps $\underline{z_\ell}(t, k), \overline{z_\ell}(t, k)$ are known for which\footnote{Here and throughout this section these inequalities are intended to hold component-wise.}  $\underline{z_\ell}(t, k)\le z_\ell(t, k)\le \overline{z_\ell}(t, k)$,
for any $(t,k)$ in the domain of the solution,
then one could define a continuous function $\lambda_\ell$ as follows
\begin{equation}
\lambda_\ell(t, k)  =  \max_{\underline{z_\ell}(t, k)\le z_\ell(t, k)\le \overline{z_\ell}(t, k)} \left\|\nabla \psi_{\ell}(z_\ell(t, k))\right\|^{2}.
\label{eq-lambda}
\end{equation}
\begin{rem} The on-line computation of (\ref{eq-lambda}) may be demanding. It may be possible to derive a simpler expression for $\lambda_\ell$ on a case-by-case basis. Suppose $n_{v_i}=1$  for any $i\in\mathcal{I}$ for instance. We can select the functions $\psi_{\ell}$ such that $\nabla \psi_{\ell}$ is nonincreasing on $\Rlo$ (take sigmoid functions for instance), (\ref{eq-lambda}) becomes then
$\lambda_{\ell}(t,k) =    \left(\nabla \psi_{\ell}(\underline{z}_\ell(t, k))\right)^{2}$ when $\underline{z}_\ell(t,k)>0$,
$\lambda_{\ell}(t,k) = \left(\nabla \psi_{\ell}(\overline{z}_\ell(t, k))\right)^{2}$ when $\overline{z}_\ell(t,k)<0$, and  $\lambda_{\ell}(t,k) = \left(\nabla \psi_{\ell}(0)\right)^{2}$ when  $\underline{z}_\ell(t,k)\overline{z}_\ell(t,k)\leq 0$. \hfill $\Box$
\end{rem}

Due to the nonlinear and distributed nature of the system, it is not an easy task to find two bounding functions $\underline{z_\ell}(t, k), \overline{z_\ell}(t, k)$ for $z_\ell(t, k)$,  unless one introduces a few additional assumptions.
\begin{ass}\label{ass-stc} The following hold.
\begin{enumerate}
\item[(i)] There exists $\bar \psi\in\Rl{}{}$  such that for any $\ell\in\{1,\ldots,M\}$ and $z_\ell\in\Rl{n_p}{}$, $\|\psi_\ell(z_\ell)\| \leq \bar \psi$.
\item[(ii)] For all $(i,j)\in \mathcal{I}^{2}$, $h_i=\mathbb{I}$ and $f_i = f_j$.
\item[(iii)] For any $i\in\mathcal{I}$,  there exist  a continuously differentiable function $V_i:\Rl{2n_{v_i}}{}$,  $\underline{\alpha}_{V_i},\overline{\alpha}_{V_i},\gamma_i\in\Kinf$ such that, for any $v_i,v_i'\in\Rl{n_{v_i}}{}$ and $u_i, u_i'\in\Rl{n_p}{}$,
\be\label{strict.incremental.passivity}
\left\{\ba{l}
\underline{\alpha}_{V_i}(\|v_i-v_i'\|)\le  V_i(v_i,v_i') \le
\overline\alpha_{V_i}(\|v_i-v_i'\|)\\[4mm]
\dst\frac{\partial V_i(v_i,v_i')}{\partial v_i}f_i(v_{i},u_{i})+\frac{\partial V_i(v_i,v_i')}{\partial v_i'}f_i(v_{i}',u_{i}')\le\\
      \hspace{2cm} -c_i V_i(v_i, v_i')+ \gamma_i (\|u_i-u_i'\|).
\ea\right.
\ee \hfill $\Box$
\end{enumerate}
\end{ass}
%
%
%
Item (i) of Assumption \ref{ass-stc} introduces no conservatism as it can be ensured by design. For example, a map $\psi_\ell$ with all the entries given by the $\arctan$ function satisfies this condition (see Section \ref{sect-simus}). The first equality of item (ii) of Assumption \ref{ass-stc} is verified by many applications, such as mechanical systems for instance where $p_i$ typically represent the position and $v_i$ the velocity, $i\in\mathcal{I}$. The second inequality simply means that the agents' dynamics are identical. The incremental input-to-state property (\cite{angeli.tac02}) in item (iii) of Assumption \ref{ass-stc} (or related concepts) is known to play a fundamental role in many problems of agreement and cooperation in dynamical networks (see \emph{e.g.}, \cite{scardovi.et.al.tac2010, stan.sepulchre.tac07}). There are interesting classes of systems for which both \eqref{eq-ass-strict-passivity} and \eqref{strict.incremental.passivity}  hold (\cite{scardovi.et.al.tac2010}); an example is provided below.

\begin{ex}\label{ex-stc}
Consider the systems of the form $\dot p_i = v_i$ and $\dot v_i = f_i(v_i,u_i)=\varphi(v_i)+u_i$, $i\in\mathcal{I}$, with the vector fields $-\varphi$ satisfying the strong monotonicity assumption\footnote{Vector fields that satisfy this property are referred to as QUAD vector fields in the literature on synchronization. The link with strict incremental passivity -- relaxed cocoercivity -- has been discussed  in \emph{e.g.}, \cite{scardovi.et.al.tac2010}. }
\begin{equation}\label{eq-ex-stc}
\begin{array}{r}
(v_i-v_i')^{\mathrm{T}} (-\varphi(v_i)+\varphi(v_i')) \ge c_i (v_i-v_i')^{\mathrm{T}} (v_i-v_i'), \quad\quad\quad \\ \forall\,v_i, v_i'\in\Rl{n_{v_i}}{},
\end{array}
\end{equation}
for some $c_i\in\Rl{}{}$. Then the storage function $S_i(v_i)=\frac{1}{2}v_i^{\mathrm{T}} v_i$ satisfies \eqref{eq-ass-strict-passivity} with $\underline{\alpha}_{S_i}(s)=\overline{\alpha}_{S_i}(s)=\frac{1}{2} s^{2}$ for any $s\in\Rlo$, and $\rho_i(v_i)=c_i \|v_i\|^{2}$ for $v_i\in\Rl{n_{v_i}}{}$, provided that $c_i\in\Rlp$ and noting that $\varphi(0)=0$ (so that $f_i(0,u_i)=0$ implies $u_i$ as required in Section \ref{sect-pb-statement}). On the other hand, $V_i(v_i, v_i')=\frac{1}{2}(v_i-v_i')^{\mathrm{T}} (v_i-v_i')$ satisfies, for any $v_i,v_i'\in\Rl{n_{v_i}}{}$ and $u_i,u_i'\in\Rl{n_p}{}$,
\begin{equation}
\begin{array}{l}
\dst\frac{\partial V_i(v_i,v_i')}{\partial v_i}(\varphi(v_{i}) + u_{i})+\frac{\partial V_i(v_i,v_i')}{\partial v_i'}(\varphi(v_{i}') + u_{i}')
    \\ \quad \quad \leq  -2c_i V_i(v_i, v_i')+ (v_i-v_i')^{\mathrm{T}} (u_i-u_i')
    \\ \quad \quad \leq  -c_i V_i(v_i, v_i')+ \frac{1}{2c_i} \|u_i-u_i'\|^{2},
\end{array}
\end{equation}
that is \eqref{strict.incremental.passivity} with $\underline{\alpha}_{V_i}(s)=\overline{\alpha}_{V_i}(s)=\frac{1}{2}s^{2}$ and $\gamma_i (s)= \frac{1}{2c_i}s^2$ for $s\in\Rlo$. \hfill $\Box$
\end{ex}


Consider the agents $i$ and $j$  connected by the edge $\ell\in\{1,\ldots,M\}$. Let $z_\ell=d_{i\ell}\, p_i + d_{j\ell}\, p_j$, $q=(p,v,\hat z,\phi)$ be a solution to (\ref{eq-sys-etc}) and $(t_{k}^{\ell},k)\in\dom q$ be such that $\phi_{\ell}(t_{k}^{\ell},k)=b_{\ell}$. We assume that no other edge triggers an event until $(t_{k+1}^{\ell},k)\in\dom q$.  We make this assumption without loss of generality only to simplify the presentation. For almost all $(t,k)\in\dom q$ with $t\geq t_{k}^{\ell}$, in view of Assumption \ref{ass-stc}, $\dot z_\ell(t, k) =\Delta v_\ell(t, k)$, where $\Delta v_\ell := d_{i\ell}\, v_i + d_{j\ell}\, v_j$.
To bound $z_\ell(t, k)$ one needs to estimate the evolution of $\Delta v_\ell(t, k)$. To this purpose,
in view of  \eqref{strict.incremental.passivity} and since $u_i=\sum_{\ell'\in \{1,\ldots,M\}}d_{i\ell}\psi_{\ell'}(\hat{z}_\ell')$, for $(t,k)\in\dom q$ with $t\geq t_{k}^{\ell}$,
{\setlength\arraycolsep{1.5pt}\begin{equation}
\begin{array}{l}
    V_i(v_i(t, k), v_j(t, k)) \leq
    \exp(-c_i(t-t_k^{\ell})) V_i (v_i(t_k^{\ell}, k), v_j(t_k^{\ell}, k))
    \\ \quad \quad  +\dst\int_{t_k^{\ell}}^{t}  \exp(-c_i(s-t_k^{\ell}))\cdot
 \\ \quad \quad \cdot\gamma_i
(
\| \dst\sum_{\ell'\in \{1,\ldots,M\}} (d_{i\ell'}-d_{j\ell'})\psi_{\ell'}(\hat{z}_{\ell'}(s,k))\|
) ds.
\end{array}\label{eq-stc-bound-V}
\end{equation}}
\hspace{-0.15cm}Using item (i) of Assumption \ref{ass-stc}, $V_i(v_i(t, k), v_j(t, k)) \leq \exp(-c_i(t-t_k^{\ell})) V_i (v_i(t_k^{\ell}, k), v_j(t_k^{\ell}, k))  +\int_{t_k^{\ell}}^{t}  \exp(-c_i(s-t_k^{\ell}))\cdot\gamma_i(2(\deg_i+\deg_j)\bar \psi)ds$.
Consequently, in view of (\ref{strict.incremental.passivity}),
\begin{equation}
\begin{array}{llllll}
\|\Delta v_\ell(t, k)\| & \leq & \overline{\Delta v_\ell}(t, k),
\end{array}
\end{equation}
with
$\overline{\Delta v_\ell}(t, k) := \underline{\alpha}_{V_i}^{-1}\Big(\!\exp(-c_i(t-t_k^{\ell}))\overline{\alpha}_{V_i}(\|\Delta v_\ell(t_k^{\ell}, k)\|)
 + \dst\frac{1}{c_i}(1-\exp(-c_i(t-t_k^{\ell})))\gamma_i (2(\deg_i$ $+\deg_j)\bar \psi)\!\Big)$.
Notice that $\overline{\Delta v_\ell}$ only depends on the value of $v_i-v_j$ at the last edge event, it is therefore available to agents $i$ and $j$ between two successive events of the edge $\ell$. One can then define the bounding maps for $z_\ell(t,k)$ as follows, for any $(t,k)\in\dom q$ with $t\geq t_{k}^{\ell}$,
\be\label{bound.on.z1}
\begin{array}{l}
\underline z_\ell(t,k):=z_\ell(t_k^{\ell},k) -\mathbf{1}\dst\int_{t_k^{\ell}}^t \overline{\Delta v_\ell}(s, k)ds \le z_\ell(t,k)
\\ \quad \quad \quad \le z_\ell(t_k^{\ell},k) +
\mathbf{1}\dst\int_{t_k^{\ell}}^t \overline{\Delta v_\ell}(s, k)ds=: \overline z_\ell(t,k).
\end{array}
\ee

\begin{rem}\label{rem-stc-relaxation-iii} The developments above indicate that item (iii) of Assumption \ref{ass-stc} can be relaxed. Indeed, the last inequality in (\ref{strict.incremental.passivity}) can be replaced by $\frac{\partial V_i(v_i,v_i')}{\partial v_i}f_i(v_{i},u_{i})+\frac{\partial V_i(v_i,v_i')}{\partial v_i'}f_i(v_{i}',u_{i}')\le \chi_i(V_i(v_i, v_i'), \|u_i-u_i'\|)$ for instance, where $\chi_i:\Rlo^{2}\to\Rl{}{}$ is non-decreasing in its last argument. In this case, we obtain, instead of (\ref{eq-stc-bound-V}), $V_i(v_i(t, k), v_j(t, k)) \leq  \eta_i(t,k)$ where $\eta_i(\cdot,k)$ is the solution to $\dot \eta_i = \chi_i(\eta_i,2(\deg_i+\deg_j) \bar\psi)$ and $\eta_i(t_k^{\ell},k)=V_i(v_i(t_k^{\ell},k),v_j(t_k^{\ell},k))$, in view of the comparison principle (see Lemma 3.4 in \cite{Khalil_02}). Note that $\eta_i$ can be computed by the agents between two successive edge events as we only need to know the values $v_i$ and $v_j$ at the last edge event to build it. We take $\overline{\Delta v_\ell}(t, k) = \underline{\alpha}_{V_i}(\eta_i(t,k))$ and follow the reasoning above.
\hfill $\Box$
\end{rem}

\subsection{Implementation of the self-triggering rules}\label{subs.st.rule}
At each event of edge  $\ell\in\{1,\ldots,M\}$, the corresponding control unit acquires the measurement $z_\ell$, it computes the control term $-d_{i\ell}\psi_\ell(z_\ell)$ as well as the next event associated with edge $\ell$. To the latter end, the control unit must compute the bounding functions $\underline{z_\ell}, \overline{z_\ell}$ according to (\ref{bound.on.z1}), the estimate $\lambda_\ell$ as in (\ref{eq-lambda}) and then solve
\begin{equation}
\begin{array}{lllll}
\dot{\phi}_\ell & = & \dst-\frac{1}{\sigma_{\ell}}\left(1+\phi_\ell^{2}\lambda_\ell\right),
\end{array}\label{eq-phi-stc}
\end{equation}
to compute the time at which $\phi_\ell$ is equal to $a_\ell$.

\subsection{Hybrid model \& analytical guarantees}

To finalize our analysis, we model the closed-loop system under self-triggering control updates  as
\begin{equation}
\begin{array}{rllll}
\left.\begin{array}{lllll}
\dot{p} & = & v\\
\dot{v} & = & f(v, \hat u)\\
\dot{\hat{z}} &=& \mathbf{0}\\
\dot{\vartheta} & = & -\mathbf{1}
\end{array}\right\} &
 \forall \ell\in\{1,\ldots,M\}\,\,\,\,\vartheta_{\ell}\geq 0\\
\left.
\begin{array}{rllll}
p^{+} & = & p\\
v^{+} & = & v\\
\left(\begin{array}{l}
\hat{z}^+\\
\vartheta^{+}
\end{array}\right) & \in &
H(z,\hat{z},\vartheta)
\end{array}\right\} & \exists \ell\in\{1,\ldots,M\}\,\,\,\,\vartheta_{\ell}=0,
\end{array}\label{eq-sys-stc}
\end{equation}
where $\vartheta:=(\vartheta_{1},\ldots,\vartheta_{M})$ and $\vartheta_{\ell}$ is a clock used to trigger the events of edge $\ell\in\{1,\ldots,M\}$. The jump map $H$ is defined similarly to (\ref{eq-jump-map-revised}) and (\ref{eq-jump-map-revised-Gamma}) {\setlength\arraycolsep{1.5pt}\begin{equation}
\begin{array}{llllll}
H(z,\hat{z},\vartheta) & := & \left\{H_{\ell}(z,\hat{z},\vartheta)\,:\,\ell\in\{1,\ldots,M\}\,\text{ and }\,\vartheta_{\ell}=0\right\}\!,
\end{array}\label{eq-jump-map-revised-stc-Gamma}
\end{equation}}
\hspace{-0.15cm}with $H_{\ell}(z,\hat{z},\vartheta) := \big(\hat{z}_{1},\ldots,\hat{z}_{\ell-1},z_{\ell},\hat{z}_{\ell+1},\ldots,\hat{z}_{M},
\vartheta_{1},\ldots,$ $ \vartheta_{\ell-1},\mathcal{T}_{\ell}(z_{\ell},\Delta v_\ell),\vartheta_{\ell+1},\ldots,\vartheta_{M}\big)$
for $\ell\in\{1,\ldots,M\}$, where $\mathcal{T}_{\ell}(z_{\ell},\Delta v_\ell)$ is the time it takes for the solution to (\ref{eq-phi-stc}) to decrease from $b_\ell$ to $a_\ell$. This constant may be analytically computed depending on the system dynamics, which helps saving CPU resources. Otherwise, (\ref{eq-phi-stc}) is solved on-line by the agents associated with edge $\ell$.

The result below is a corollary of Theorem \ref{th-etc-N}.
\begin{cor}\label{cor-stc-N} Consider system (\ref{eq-sys-stc}) and suppose the following holds.
\begin{enumerate}
\item[(i)] Items (i)-(iii) of Theorem \ref{th-etc-N} hold.
\item[(ii)] Assumption \ref{ass-stc} is guaranteed.
\end{enumerate}
The solutions have a uniform semiglobal average dwell-time and the maximal solutions are complete and approach the set\footnote{In the definition of the set, $\vartheta_\ell\in[0,\sigma_{\ell}(b_{\ell}-a_{\ell})]$. This comes from the fact that the inter-edge event times are less than or equal the time it takes for $\dot \omega_\ell = -\frac{1}{\sigma_\ell}$ to decrease from $b_\ell$ to $a_\ell$ (in view of the comparison principle), which is equal to $\sigma_{\ell}(b_{\ell}-a_{\ell})$.} $\left\{(p,v,\hat{z},\vartheta): z\in\mathcal{A},\,v=\mathbf{0},\,\vartheta_{\ell}\in[0,\sigma_{\ell}(b_{\ell}-a_{\ell})]\right.$ for $\left.\ell\in \{1,\ldots,M\}\right\}$. \hfill $\Box$ 
\end{cor}


\section{Simulation results}\label{sect-simus}

\begin{figure}
\centering
\psfrag{1}[][][1]{\small $1$}
\psfrag{2}[][][1]{\small $2$}
\psfrag{3}[][][1]{\small $3$}
\psfrag{4}[][][1]{\small $4$}
\psfrag{5}[][][1]{\small $5$}
 \includegraphics[scale=0.25]{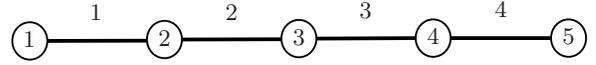}
 \vspace{-0.3cm}\caption{Graph in Section \ref{sect-simus} ($N=5$ and $M=4$).}\label{fig-graph}
\end{figure}

The objective is to ensure the rendez-vous of $N=5$ identical agents when the graph is given by a line as depicted in Figure \ref{fig-graph}. The number of edges is $M=4$ in this case. We consider the following agents' dynamics which are subject to {\em input saturation}
\begin{equation}
\begin{array}{lllll}
\dot{p}_{i} & = & v_{i}\\
\dot{v}_{i} & = & - v_i + u_i,
\end{array}\label{eq-sys-ex}
\end{equation}
where $p_i\in\Rl{}{}$, $v_i\in\Rl{}{}$, $u_i\in[-\bar u,\bar u]$ is the saturated control input and $\bar u = 1$ is the saturation level. We want to ensure the rendez-vous of the agents, in other words we want the relative distances $z_{\ell}$ to converge to the origin, hence $\mathcal{A}_\ell = \{0\}$, $\ell\in\{1,\ldots,M\}$.

System (\ref{eq-sys-ex}) verifies Assumption \ref{ass-strict-passivity} with $S_i(v_i)=\frac{1}{2} v_i^{2}$, $\underline{\alpha}_{S_i}(s)=\overline{\alpha}_{S_i}(s)=\frac{1}{2} s^{2}$, and $\rho_i(v_i)=v_i^{2}$ for $i\in\{1,\ldots,N\}$, $v_i\in\Rl{}{}$ and $s\in\Rlo$. We design the control input $u_i$ as in (\ref{eq-cont}) with $\psi_{\ell}(z_\ell)=\frac{1}{\pi} \arctan (z_{\ell})$ for $z_{\ell}\in\Rl{}{}$ and $\ell\in\{1,\ldots,M\}$. Hence $P_{\ell}(z_{\ell})=\frac{1}{\pi}\left(z_{\ell}\arctan(z_{\ell})-\frac{1}{2} \ln(1+z_{\ell}^{2})\right)$. We see that items (a), (b) and (d) in Section \ref{sect-pb-statement} are verified. Noting that $P_{\ell}$ is positive definite, continuous and radially unbounded, we apply Lemma 4.3 in \cite{Khalil_02} to deduce that item (c) in Section \ref{sect-pb-statement} holds. We notice that this choice of $\psi_{\ell}$, $\ell\in\{1,\ldots,M\}$, ensures that all the control inputs lie in the admissible range $[-1,1]$ as $\frac{1}{\pi}\arctan(\Rl{}{})=(-\frac{1}{2},\frac{1}{2})$ and the maximal degree of the agents is $2$ (see Figure \ref{fig-graph}).

Our aim is to design event-triggered, time-triggered and self-triggered controllers. We first concentrate on the synthesis of the event-triggered controllers. We therefore need to verify that the conditions of Theorem \ref{th-etc-N} hold. 
We select $\sigma_{\ell}=\frac{\kappa_i}{4}$ with $\kappa_i=\frac{1}{4}$, $a_{\ell}=0$, different values will be assigned to $b_{\ell}$, and we initialize the clock variables $\phi_{ij}$ at the same values, so that Assumption \ref{ass-synchro} a fortiori holds. Hence item (i) of Theorem \ref{th-etc-N} is ensured. Noting that in our case $h_{i}(v_i)=v_i$ for $v_i\in\Rl{}{}$, we have $\|h_i(v_i)\|^{2}=\rho_i(v_i)$  and (\ref{eq-thm-ass-cond}) holds. Our choice of $\sigma_{\ell}$, $\ell\in\{1,\ldots,M\}$, guarantees (\ref{eq-thm-ass-param}), as consequence item (ii) of Theorem \ref{th-etc-N} is ensured. We note that item (iii) of Theorem \ref{th-etc-N} applies since $\mathcal{A}=\{0\}^{5}$ (see Section III in \cite{Arcak-tac07}). Consequently, the conclusions of Theorem \ref{th-etc-N} hold.
To design time-triggered controllers, we also need to ensure (\ref{eq-boundedness-nabla-psi}), which is the case by taking $K_{\ell}=\frac{1}{\pi}$ for $\ell\in\{1,\ldots,M\}$. We have selected $T_\ell$ as in (\ref{eq-T}) and $\varepsilon_{\ell}=T_{\ell}$, which means that each sequence of edge events is $T_{\ell}$-periodic. Finally, we verify that Assumption \ref{ass-stc} is verified by system (\ref{eq-sys-ex}) for the construction of the self-triggered controllers. Items (i)-(ii) of Assumption \ref{ass-stc} hold with $\bar \psi = \frac{1}{2}$. Item (iii) of Assumption \ref{ass-stc} is verified in view of Example \ref{ex-stc} as (\ref{eq-ex-stc}) holds with $c_i=1$.

\begin{figure}[htbp]
 \begin{center}
 \includegraphics[scale=0.25]{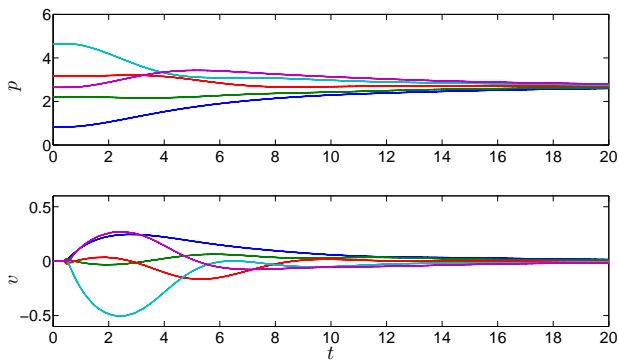}
 \end{center}
\caption{Evolution of $p_i$ (top) and $v_i$ (bottom) for $b=10$ using the event-triggered controllers, $i\in\{1,\ldots,N\}$.\vspace{0cm}}
\label{fig-ex}
\end{figure}

An example of the evolution of $p_i$ and $v_i$, $i\in\{1,\ldots,N\}$, is provided in Figure \ref{fig-ex}, which has been obtained by using the event-triggered controllers with $b=10$. We see that the rendez-vous is ensured and that all the $v_i$'s converge to the origin as expected. We have then simulated the system with the three types of controllers for $100$ different initial conditions for which $p$ is randomly distributed in $[0,5]$, $v(0,0)=\mathbf{0}$, $\hat{z}(0,0)=D^{\mathrm{T}}p(0,0)$, $\phi(0,0)=b\mathbf{1}$, $\vartheta_\ell(0,0) = \mathcal{T}_\ell(z_\ell(0,0),\Delta v_\ell(0,0))$ and $\tau_{\ell}(0,0)$ is randomly distributed in $[0,T_{\ell}]$ for the time-triggered controllers (in order to avoid synchronous periodic events over the whole network), with a simulation time of $20$s and for different vlaues of $b$. Table \ref{tab-ex} provides the obtained averages of the total number of edge events, and the averages of the time it takes for $\|z\|=\|(z_1,\ldots,z_M)\|$ to become less than $5\%$ of its initial value, which we denote $t_{5\%}$ and which serves as a measure of the speed of convergence. The results show that the event-triggered controllers generally generate less events compared to the self-triggered controllers, however the difference is not significant, which justifies the proposed design method of the self-triggered controllers in Section \ref{sect-stc}. Also, the self-triggered controllers give rise to less events compared to the time-triggered controllers, which is in agreement with the theoretical developments. On the other hand, less events typically leads to longer times $t_{5\%}$, which can be explained by the fact that the control inputs are more often updated and the states thus converge faster. Table \ref{tab-ex} also suggests that to increase the value of $b$ reduces the number of edge events at the price of a longer convergence time. The parameter $b_{\ell}$ (equivalently $a_{\ell}$ and $\sigma_\ell$), $\ell\in\{1,\ldots,M\}$, may therefore be adjusted to reduce the communication and computation cost at the price of a slower convergence speed.

\begin{table}[t]
{\scriptsize\centering
\begin{tabular}{cc|ccc}
& & $b=1$ & $b=10$ & $b=100$ \\
\hline & & & & \vspace{-0.3cm}\\
Average \# of events & ETC & $1313.8$ & $291.29$ & $219.84$\\
                     & STC & $1313.7$ & $292.58$ & $224.35$\\
                     & TTC & $1322.1$ & $321.49$ & $264.60$\\
\hline & & & & \vspace{-0.3cm}\\
Average $t_{5\%}$   & ETC & $11.782$   & $13.1884$ & $15.4087$\\
                     & STC & $11.924$  & $12.8762$ & $13.6525$\\
                     & TTC & $13.0180$ & $11.7173$ & $12.3144$\\
\end{tabular}
\caption{Average number of events and average value of $t_{5\%}$ (\#: number, ETC: event-triggered control, STC: self-triggered control, TTC: time-triggered control).\vspace{0cm}}
\label{tab-ex}}
\end{table}

\section{Proof of Theorem \ref{th-etc-N}} \label{sect-proof-thm-etc}

For the sake of convenience, we write system (\ref{eq-sys-etc}) as
\begin{equation}
\begin{array}{llllllllllllllll}
\dot{q} = F(q) & \text{ for } q\in C, & &
q^{+} \in G(q) & \text{ for }q\in D,
\end{array}\label{eq-sys-etc-q}
\end{equation}
where $q:=(p,v,\hat{z},\phi)\in\Rl{n_q}{}$,
$C := \big\{q \,:\, \forall \ell\in \{1,\ldots,M\}\,\,\,\, \phi_{\ell}\in[a_{\ell},b_{\ell}]\big\}$,  $D := \big\{q \,:\, \exists \ell\in\{1,\ldots,M\} \,\,\,\, \phi_{\ell}=a_{\ell}\big\}$, and $n_q:=n_p N + n_v + n_p M +M$.

\subsection{Lyapunov analysis}\label{subsect-proof-lyap}

The analysis is performed relying on Lyapunov arguments. To this purpose, we introduce the function
\begin{equation}
\begin{array}{llllll}
U(q):= U_{\text{phys}}(q)+U_{\text{cyber}}(q) & & \forall q\in \Rl{n_q}{}.
\end{array}\label{lyap}
\end{equation}
The term $U_{\text{phys}}$ takes into account the physical component of the system and is an energy-like function of the form
\begin{equation}
\begin{array}{llllll}
U_{\text{phys}} (q):= \dst  \sum_{i\in \mathcal{I}} S(v_i) + \sum_{\ell\in \{1,\ldots,M\}} P_\ell(z_{\ell})  & & \forall q\in \Rl{n_q}{},
\end{array}\label{lyap.phys}
\end{equation}
where $S_i$ and $P_{\ell}$ respectively come from Assumption \ref{ass-strict-passivity} and the definition of $\psi_{\ell}$ in Section \ref{sect-pb-statement}.
The term $U_{\text{cyber}}$  takes into account the cyber-physical nature of the system  and it will be specified in the following. For $q\in C$, one obtains from Assumption \ref{ass-strict-passivity}
\begin{equation}
\begin{array}{llllll}
\langle \nabla U_{\text{phys}}(q), F(q)\rangle  & = &  \dst\sum_{i\in \mathcal{I}} \left\{-\rho_i(v_i)+\hat
u_i^{\mathrm{T}}y_i\right\}
\\
&& +  \sum_{\ell\in \{1,\ldots,M\}} \left\langle \nabla P_{\ell}(z_{\ell}), \Delta y_{\ell}\right\rangle
\end{array}
\end{equation}
where $\Delta y_\ell:=y_{j}-y_{i}$ when $j$ positive end of the edge $\ell$ and
$\Delta y_\ell:=y_{i}-y_{j}$ when $i$ positive end of the edge $\ell$,
and where we have exploited the fact that $d_{i\ell} y_i+d_{j\ell} y_j=\Delta y_\ell$.
By definition of $\psi_{\ell}$,
\begin{equation}
\begin{array}{llllll}
\langle \nabla U_{\text{phys}}(q), F(q)\rangle  & = &  \dst\sum_{i\in \mathcal{I}}
\left\{-\rho_i(v_i)\hat u_i^{\mathrm{T}}y_i\right\}
\\ && + \sum_{\ell\in \{1,\ldots,M\}} \left\langle \psi_{\ell}(z_{\ell}), \Delta y_{\ell}\right\rangle.
\end{array}
\end{equation}
We write $\dst\sum_{\ell\in \{1,\ldots,M\}} \left\langle \psi_{\ell}(z_{\ell}), \Delta y_{\ell}\right\rangle=\Psi(z)^{\mathrm{T}}\Delta y$ where $\Psi(z)=\left(\psi_{1}(z_1),\ldots,\psi_{M}(z_M)\right)$, and $\Delta y:=(\Delta y_1,\ldots,\Delta y_M)$. We have  $\Delta y = (D^{\mathrm{T}}\otimes \mathbb{I}_{n_{p}})y$ where $y:=(y_1,\ldots,y_M)$. Therefore
$\dst\sum_{\ell\in \{1,\ldots,M\}} \left\langle \psi_{\ell}(z_{\ell}),\right.$ $ \left. \Delta y_{\ell}\right\rangle  = \Psi(z)^{\mathrm{T}}(D^{\mathrm{T}}\otimes \mathbb{I}_{n_{p}})y
= \Psi(z)^{\mathrm{T}}(D \otimes \mathbb{I}_{n_{p}})^{\mathrm{T}}y$,
where $(D^{\mathrm{T}}\otimes \mathbb{I}_{n_{p}})=(D \otimes \mathbb{I}_{n_{p}})^{\mathrm{T}}$ is used to obtain the last equality. Noticing that $u :=(u_1,\ldots,u_N)=-(D\otimes \mathbb{I}) \Psi(z)$, we obtain $\dst\sum_{\ell\in \{1,\ldots,M\}} \left\langle \psi_{\ell}(z_{\ell}), \Delta y_{\ell}\right\rangle = -u^{\mathrm{T}}y$. As a consequence,
\begin{equation}
\begin{array}{llllll}
\langle \nabla U_{\text{phys}}(q), F(q)\rangle  & = &  \dst\sum_{i\in \mathcal{I}} \left\{-\rho_i(v_i)+\hat u_i^{\mathrm{T}}y_i\right\}
-u^{\mathrm{T}}y\\
&=&  \dst\sum_{i\in \mathcal{I}} \left\{-\rho_i(v_i)+(\hat u_i-u_i)^{\mathrm{T}}y_i\right\}.
\end{array}\label{phys.deriv}
\end{equation}
The interpretation of the expression (\ref{phys.deriv}) is clear. The use of sampled-data measurements $\hat z_\ell$, instead of the actual measurements $z_\ell$, causes the appearance of a perturbative term $\dst\sum_{i\in \mathcal{I}}(\hat u_i-u_i)^{\mathrm{T}}y_i$ in the derivative of the energy function $U_{\text{phys}}$, potentially disrupting the achievement of the coordination. How this perturbation can be counteracted is explained by the introduction of  the term $U_{\text{cyber}}$ in the Lyapunov function in (\ref{lyap}) 
\begin{equation}
U_{\text{cyber}}(q):= \dst\sum_{\ell\in  \{1,\ldots,M\}} \frac{1}{2}\phi_\ell  \left\|\psi_\ell(\hat z_\ell)-\psi_\ell(z_\ell)\right\|^{2} \quad \forall q\in \Rl{n_q}{}.
\label{lyap.cyber}
\end{equation}
We show that the update law for  $\phi_\ell$ guarantees that the Lyapunov function $U$ does not increase as far as $q\in C$. In fact, observe that, for $q\in C$,
\begin{equation}
\begin{array}{llllll}
\langle \nabla U_{\text{cyber}}(q), F(q)\rangle =
   \\ \!\dst\sum_{\ell\in \{1,\ldots,M\}} \!\!\!\Big\{-\frac{1}{2\sigma_{\ell}}\left(1+\phi_{\ell}^{2}\left\|\nabla \psi_{\ell}(z_{\ell})\right\|^{2}\right)  \left\|\psi_\ell(\hat z_\ell)-\psi_\ell(z_\ell)\right\|^2
   \\  -\phi_\ell (\psi_\ell(\hat z_\ell)-\psi_\ell(z_\ell))^{\mathrm{T}} \nabla \psi_\ell(z_\ell) \Delta y_\ell\Big\}.
\end{array}
\end{equation}
The last term on the right-hand side above is upper-bounded as follows
$-\sum_{\ell\in \{1,\ldots,M\}}\phi_\ell (\psi_\ell(\hat z_\ell)-\psi_\ell(z_\ell))^{\mathrm{T}} \nabla \psi_\ell(z_\ell) \Delta y_\ell \leq
\sum_{\ell\in \{1,\ldots,M\}}\Big\{\frac{1}{2\sigma_{\ell}}\phi_\ell^{2}\left\|\psi_\ell(\hat z_\ell)-\psi_\ell(z_\ell)\right\|^{2} \left\|\nabla \psi_\ell(z_\ell)\right\|^{2} + \frac{\sigma_{\ell}}{2}\left\|\Delta y_\ell\right\|^{2}\Big\}$.
Consequently,
\begin{equation}
\begin{array}{llllll}\label{cyber.deriv}
\langle \nabla U_{\text{cyber}}(q), F(q)\rangle  = \\[4mm]
 \dst\sum_{\ell\in \{1,\ldots,M\}} \Big\{-\frac{1}{2\sigma_{\ell}} \left\|\psi_\ell(\hat z_\ell)-\psi_\ell(z_\ell)\right\|^2 + \frac{\sigma_{\ell}}{2}\left\|\Delta y_\ell\right\|^{2}\Big\}.
\end{array}
\end{equation}
Overall, from (\ref{phys.deriv}) and (\ref{cyber.deriv}),
\begin{equation}
\ba{lllll}
\langle \nabla U(q), F(q)\rangle  \leq   \dst\sum_{i\in \mathcal{I}} \left\{-\rho_i(v_i)+(\hat u_i-u_i)^{\mathrm{T}}y_i\right\} \\
 + \dst\sum_{\ell\in \{1,\ldots,M\}} \Big\{-\frac{1}{2\sigma_{\ell}} \left\|\psi_\ell(\hat z_\ell)-\psi_\ell(z_\ell)\right\|^2 + \frac{\sigma_{\ell}}{2}\left\|\Delta y_\ell\right\|^{2}\Big\}.
\ea\end{equation}
Since $u=-(D\otimes \mathbb{I})\Psi(z)$ and $\hat u=-(D\otimes \mathbb{I})\Psi(\hat z)$, we obtain
$\dst\sum_{i\in \mathcal{I}} (\hat u_i-u_i)^{\mathrm{T}}y_i=(\hat u - u)^{\mathrm{T}}y = -(\Psi(\hat z)-\Psi(z))^{\mathrm{T}}(D\otimes \mathbb{I})^{\mathrm{T}}y=-(\Psi(\hat z)-\Psi(z))^{\mathrm{T}}(D^{\mathrm{T}}\otimes \mathbb{I})y=-(\Psi(\hat z)-\Psi(z))^{\mathrm{T}}\Delta y$ (as $\Delta y = (D^{\mathrm{T}}\otimes \mathbb{I}_{n_{p}})y$). Consequently,
\begin{equation}
\label{total.deriv-prel}
\ba{lllll}
\langle \nabla U(q), F(q)\rangle  \leq  \\ \dst -\sum_{i\in \mathcal{I}} \rho_i(v_i) - \sum_{\ell\in\{1,\ldots,M\}}(\psi_{\ell}(\hat z_{\ell})-\psi_{\ell}(z_{\ell}))^{\mathrm{T}}\Delta y_{\ell}\\
 + \dst\sum_{\ell\in \{1,\ldots,M\}} \Big\{-\frac{1}{2\sigma_{\ell}} \left\|\psi_\ell(\hat z_\ell)-\psi_\ell(z_\ell)\right\|^2 + \frac{\sigma_{\ell}}{2}\left\|\Delta y_\ell\right\|^{2}\Big\}.
\ea
\end{equation}
We use the inequality
$\dst -\sum_{\ell\in\{1,\ldots,M\}}(\psi_{\ell}(\hat z_{\ell})-\psi_{\ell}(z_{\ell}))^{\mathrm{T}}\Delta y_{\ell} \leq \dst\sum_{\ell\in\{1,\ldots,M\}}\Big\{\frac{1}{2\sigma_{\ell}} \left\|\psi_{\ell}(\hat z_{\ell})-\psi_{\ell}(z_{\ell})\right\|^{2} + \frac{\sigma_{\ell}}{2}\left\|\Delta y_{\ell}\right\|^{2}\Big\}
$ to obtain from (\ref{total.deriv-prel})
\begin{equation}
\ba{lllll}
\langle \nabla U(q), F(q)\rangle & \leq  & \dst -\sum_{i\in \mathcal{I}} \rho_i(v_i) + \sum_{\ell\in\{1,\ldots,M\}}\sigma_{\ell}\left\|\Delta y_{\ell}\right\|^{2}.
\ea
\end{equation}
Bearing in mind that $\|\Delta y_\ell\|^2 \le 2\|y_i\|^2 +2\|y_j\|^2$, the latter term satisfies
\begin{equation}
\begin{array}{lllll}
\dst\sum_{\ell\in \{1,\ldots,M\}} \sigma_{\ell} \|\Delta y_\ell\|^2 & \le & \dst\sum_{i\in \mathcal{I}} 2 \deg_i\max_{\ell\in \mathcal{E}_i}\{\sigma_\ell\}\|y_i\|^2,
\end{array}
\end{equation}
and we have
\begin{equation}
\langle \nabla U(q), F(q)\rangle \leq   \dst -\sum_{i\in \mathcal{I}} \rho_i(v_i) + 2\sum_{i\in\mathcal{I}}\deg_i \max_{\ell\in\mathcal{E}_i}\{\sigma_{\ell}\}\left\|y_{i}\right\|^{2}.
\end{equation}
We finally use (\ref{eq-thm-ass-cond}) to derive
\begin{equation}
\label{eq-proof-thm-U-flow}
\ba{lllll}
\langle \nabla U(q), F(q)\rangle  \leq \\
   \quad \quad \quad \dst -\sum_{i\in \mathcal{I}} (1-\kappa_i)\rho_i(v_i) \leq -(1-\kappa)\sum_{i\in \mathcal{I}}\rho_i(v_i),
\ea
\end{equation}
where $\kappa:=\dst\max_{i\in\mathcal{I}}\kappa_i>0$.

Let $q\in D$, $U_{\text{phys}}(G(q)) = U_{\text{phys}}(q)$ since the function $U_{\text{phys}}$ only includes  $z$  and $v$ that do not undergo jumps. On the other hand, the term $U_{\text{cyber}}$ satisfies, in view of (\ref{eq-sys-etc}),
\begin{equation}
\begin{array}{lllllll}
U_{\text{cyber}}(G(q))  =  \dst\sum_{\ell\in  \{1,\ldots,M\}\backslash\{\ell'\}} \frac{1}{2}\phi_\ell  \left\|\psi_\ell(\hat{z}_\ell)-\psi_\ell(z_\ell)\right\|^2
  \\ \quad\quad\quad \dst +\frac{1}{2} b_{\ell'}  \left\|\psi_\ell(z_{\ell'})-\psi_\ell(z_{\ell'})\right\|^2,
\end{array}
\end{equation}
where $\ell'\in\{1,\ldots,M\}$ is such that $\phi_{\ell'}=a_{\ell'}$. As a consequence
$U_{\text{cyber}}(G(q)) = \dst\sum_{\ell\in  \{1,\ldots,M\}\backslash\{\ell'\}} \frac{1}{2} \phi_\ell $ $\times\left\|\psi_\ell(\hat{z}_\ell)-\psi_\ell(z_\ell)\right\|^2 \leq U_{\text{cyber}}(q)$.
Hence, we conclude that
\begin{equation}
\begin{array}{lllllll}\label{eq-proof-thm-U-jump}
U(G(q)) & \leq & U(q).\\
\end{array}
\end{equation}

\subsection{Completeness and boundedness properties of the maximal solutions}\label{subsect-proof-max-complete}

We now use the conclusions of Section \ref{subsect-proof-lyap} to prove the completeness of the maximal solutions to (\ref{eq-sys-etc-q}) as well as some boundedness properties which will be essential in the sequel.

We first show that any maximal solution to (\ref{eq-sys-etc-q}) is nontrivial. We verify for that purpose that $F(q)\in T_{C}(q)$ for any $q\in C\backslash D$ in view of Proposition 6.10 in \cite{Goebel-Sanfelice-Teel-book}, where $T_{C}(q)$ is the tangent cone to $C$ at $q$ (see Definition \ref{def-tangent-cone}). Let $q\in C\backslash D$, if $q$ is the interior of $C$, $T_{C}(q)=\Rl{n_q}{}$ and the desired property holds. If $q\in C\backslash D$ and $q$ is not in the interior of $C$, then necessarily there exists $\ell\in\{1,\ldots,M\}$ such that $\phi_{\ell}=b_{\ell}$. We suppose that there is a unique such $\ell$ for the sake of simplicity (similar arguments apply when it is not the case). In this case, $T_{C}(q)=\Rl{n_q-M}{}\times\Rl{\ell-1}{}\times(-\infty,0]\times\Rl{M-\ell}{}$ and $F(q)\in T_{C}(q)$ as the flow map of $\phi_{\ell}$ at $q$ is strictly negative in view of (\ref{eq-phi-l}). Consequently, any maximal solution to (\ref{eq-sys-etc-q}) is nontrivial.

Let $q$ be a maximal solution to (\ref{eq-sys-etc-q}). Since $G(D)\subset C$, we know from Proposition 6.10 in \cite{Goebel-Sanfelice-Teel-book} that we only need to prove that $q$ does not explode in finite (hybrid) time to ensure that $q$ is complete. As a consequence of Assumption \ref{ass-strict-passivity} and item (c) in Section \ref{sect-pb-statement},
for any $q\in\Rl{n_q}{}$,
$\dst\sum_{i\in\mathcal{I}}\underline{\alpha}_{S_i}(\|v_i\|) + \sum_{\ell\in\{1,\ldots,M\}}\underline{\alpha}_{P_\ell}(\|z_\ell\|_{\mathcal{A}_{\ell}}) \leq U(q)$.
Noting that $\|z\|_{\mathcal{A}} \leq \dst\sum_{\ell\in\{1,\ldots,M\}}\|z_\ell\|_{\mathcal{A}_{\ell}}$ for any $z=(z_1,\ldots,z_M)$, and using Remark 2.3 in \cite{Laila-Nesic-cdc-02(small-gain)}, we deduce that there exists $\underline{\alpha}_{U}\in\Kinf$ such that $\underline{\alpha}_{U}(\|(z,v)\|_{\mathcal{A}\times\{0\}^{n_v}}) \leq U(q)$ for any $q\in\Rl{n_q}{}$.
We know that $U$ does not increase along the solutions to (\ref{eq-sys-etc-q}) in view of (\ref{eq-proof-thm-U-flow}) and (\ref{eq-proof-thm-U-jump}), thus, for all $(t,k)\in\dom q$,
\begin{equation}
\begin{array}{rllllll}
\|(z(t,k),v(t,k))\|_{\mathcal{A}\times\{0\}^{n_v}} & \leq & \underline{\alpha}_{U}^{-1}(U(q(0,0))),
\end{array}
\end{equation}
which implies that there exists a constant $\Theta(q(0,0))\in\Rlo$ such that, for any $(t,k)\in\dom q$,
\begin{equation}
\begin{array}{rllllll}
\|(z(t,k),v(t,k))\| & \leq & \Theta(q(0,0))
\end{array}\label{eq-proof-thm-bound-z-v}
\end{equation}
since $\mathcal{A}\times\{0\}^{n_v}$ is a compact set.
Consequently, in view of (\ref{eq-sys-etc-q}),
\begin{equation}
\begin{array}{rllllll}
\|\hat{z}(t,k)\| & \leq & \widehat\Theta(q(0,0))
\end{array}\label{eq-proof-thm-bound-z-hat}
\end{equation}
for some $\widehat\Theta(q(0,0))\geq 0$ and any $(t,k)\in\dom q$. Noting that $\phi(t,k)\in [a_{1},b_{1}]\times\ldots\times[a_{M},b_{M}]$ for any $(t,k)\in\dom q$, we are left with proving that $p$ does not explode in finite time. At each jump, $p$ does not vary. On flows, we have $\dot p = y$. Since $y=h(v):=(h_1(v_1),\ldots,h_N(v_N))$, $h$ is continuous (as it is locally Lipschitz) and $v$ is ensured to be bounded in view of (\ref{eq-proof-thm-bound-z-hat}), $p$ may grow at least linearly during flows, which guarantees that it does not explode in finite time. Therefore, by Proposition 6.10 in \cite{Goebel-Sanfelice-Teel-book}, we know that $q$ is complete. Note that we do not guarantee a boundedness property for $p$ contrary to the other variables: that is not needed to ensure the desired coordination objective as we will see. 

\subsection{Auxiliary system}\label{subsect-proof-aux-sys}

The invariance principle in Theorem 8.2 in \cite{Goebel-Sanfelice-Teel-book} applies to precompact solutions of the considered hybrid system, \ie to maximal solutions which are complete and for which the closure of their range is bounded. Completeness of the maximal solutions to (\ref{eq-sys-etc-q}) has been established in Section \ref{subsect-proof-max-complete}, however we have not proved the required boundedness property because of the $p$-component of the solutions. We overcome this issue by considering the auxiliary system below
\begin{equation}
\begin{array}{rllll}
\left.\begin{array}{lllll}
\dot{z} & = & (D^{\mathrm{T}}\otimes \mathbb{I})y\\
\dot{v} & = & f(v,\hat u) \\
\dot{\hat{z}} &=& \mathbf{0}\\
\dot{\phi} & = & -\Sigma^{-1}
\left(\mathbf{1}+ \Phi(z)\right)\\
\end{array}\!\!\!\!\right\} &
\hspace{-0.2cm} \forall \ell\in\{1,\ldots,M\}\,\,\,\phi_{\ell}\in [a_{\ell},b_{\ell}]\\
\left.
\begin{array}{rllll}
z^{+} & = & z\\
v^{+} & = & v\\
\left(\begin{array}{l}
\hat{z}^+\\
\phi^{+}
\end{array}\right) & \in &
G(z,\hat{z},\phi)
\end{array}\!\!\!\!\right\} & \hspace{-0.2cm}\exists \ell\in\{1,\ldots,M\}\,\,\,\phi_{\ell}=a_{\ell},
\end{array}\label{eq-sys-etc-aux}
\end{equation}
which we denote by, for the sake of convenience,
\begin{equation}
\begin{array}{llllllllllllllll}
\dot{q}_{\aux} = F_{\aux}(q_{\aux}) & \text{ for } q_{\aux}\in
C_{\aux},\\
q_{\aux}^{+} \in G_{\aux}(q_{\aux}) & \text{ for }q_{\aux}\in D_{\aux},
\end{array}\label{eq-sys-etc-q-aux}
\end{equation}
where $q_{\aux}:=(z,v,\hat{z},\phi)\in\Rl{n_{q_\aux}}{}$ and $C_{\aux} := \big\{q_{\aux} \,:\, \forall \ell \in \{1,\ldots,M\}\,\,\,\, \phi_{\ell}\in[a_{\ell},b_{\ell}]\big\}$, $D_{\aux} := \big\{q_{\aux} \,:\, \exists \ell\in\{1,\ldots,M\} \,\,\,\, \phi_{\ell}=a_{\ell}\big\}$, $n_{q_{\aux}}:=n_p M+n_v+n_p M + M$.
The difference with (\ref{eq-sys-etc}) is that the state $p$ has been replaced by the relative distance $z$, while the other state variables remain unchanged. This change of variable is not invertible: $z=(D^{\mathrm{T}}\otimes\mathbb{I})p$ and $\text{rank}(D^{\mathrm{T}}\otimes \mathbb{I})=\text{rank}(D)\times \text{rank}(\mathbb{I})=(N-1)n_p\neq N n_p $ ($\text{rank}(D)=N-1$ since the graph is connected, see Theorem 8.3.1 in \cite{Godsil-Royle-2001(graph-book)}). Nevertheless, we argue that to prove the desired convergence property on system (\ref{eq-sys-etc-aux}) ensures the same property holds for system (\ref{eq-sys-etc}). Indeed, in (\ref{eq-sys-etc-aux}), the flow and the jump maps of $p$ and of $(z,v,\hat z,\phi)$ are decoupled. We can thus isolate the dynamics of these two systems and only study the latter, provided that the maximal solutions to the $p$-system are complete, which is the case in view of Section \ref{subsect-proof-max-complete}.

We will therefore apply an hybrid invariance principle in Chapter 8 of \cite{Goebel-Sanfelice-Teel-book} to system (\ref{eq-sys-etc-aux}). We first note that this system is (nominally) well-posed for the same reasons as system (\ref{eq-sys-etc}) is. Furthermore, the maximal solutions to (\ref{eq-sys-etc-aux}) are complete in view of Section \ref{subsect-proof-max-complete} and the closure of their range is bounded in view of (\ref{eq-proof-thm-bound-z-v}), (\ref{eq-proof-thm-bound-z-hat}) and the fact that $\phi\in [a_1,b_1]\times\ldots\times[a_M,b_M]$. Thus, the maximal solutions to (\ref{eq-sys-etc-aux}) are precompact.

\subsection{Average dwell-time solutions}\label{subsect-proof-avg-dwt}
Next step is to show that the solutions to (\ref{eq-sys-etc-aux}) have a uniform semiglobal average dwell-time (see Definition \ref{def-prel-adt}). This property is important for practical reasons as explained in Section \ref{sect-pb-statement}, furthermore it will be useful to prove the desired convergence property.
To this end, we first study the time interval between two successive events associated with a given edge. In other words, we investigate the time it takes for the clock $\phi_{\ell}$ to decrease from $b_{\ell}$ to $a_{\ell}$ in view of (\ref{eq-phi-l}), for $\ell\in\{1,\ldots,M\}$.

Let $\Delta>0$ and take a solution $q_\aux$ to (\ref{eq-sys-etc-q-aux}) such that $\|q_\aux(0,0)\|\leq \Delta$. According to (\ref{eq-proof-thm-bound-z-v}), there exists $\bar{\Delta}>0$ (which depends on $\Delta$) such that, for any $(t,k)\in\dom  q_\aux$,
for any $\ell\in\{1,\ldots,M\}$,
$\|z_{\ell}(t,k)\| \leq \|(z(t,k),v(t,k))\| \leq \bar{\Delta}$.
On the other hand, the time it takes from $\phi_{\ell}$ to decrease from $b_{\ell}$ to $a_{\ell}$ is lower bounded by the time it takes for $\theta_{\ell}$, the solution to the differential equation below
\begin{equation}
\begin{array}{llllllllllll}
\dot{\theta}_{\ell} & = & -\dst\frac{1}{\sigma_{\ell}}\left(1+ \theta_{\ell}^{2} { \underset{\xi \text{ s.t.} \|\xi\| \leq \bar{\Delta}}\max}\left\|\nabla \psi_{\ell}(\xi)\right\|^{2}\right),
\end{array}\label{eq-theta}
\end{equation}
to decrease from $b_{\ell}$ to $a_{\ell}$, in view of (\ref{eq-phi-l}) and according to the comparison principle (see Lemma 3.4 in \cite{Khalil_02}). Note that the maximum in (\ref{eq-theta}) is well-defined since $\|\nabla \psi_{\ell}\|^{2}$ is continuous and since it is taken over a compact set. The aforementioned time interval\footnote{In fact, the rate of change of $\theta_{\ell}$ is upper and lower bounded as follows
\[
-\frac{1}{\sigma_{\ell}}\left(1+b_{\ell}^2  {\underset{\xi \text{ s.t. } \|\xi\|\leq \bar{\Delta}}\max}\big\{\left\|\nabla \psi_{\ell}(\xi)\right\|^{2}\}\right) \le
\dot \theta_{\ell} \le  -\frac{1}{\sigma_{\ell}}.
\]
}
is obviously a strictly positive constant $\tau_{\ell}(a_{\ell},b_{\ell},\Delta)$ in view of (\ref{eq-theta}) (recall that $\bar{\Delta}$ depends on $\Delta$). Consequently, the ordinary time between two successive events associated with the edge $\ell$ is lower bounded by $\tau_{\ell}(a_{\ell},b_{\ell},\Delta)$. Let $q_\aux$ be a solution
to  (\ref{eq-sys-etc-aux}) and $(s,i),(t,k)\in\dom q_\aux$ with $s+i\leq t+k$. In view of the above developments, the number of events associated with  the edge $\ell$ between $(s,i)$ and $(t,k)$, which can be written as a function $\mathbf{n}_{\ell}(s,t)$ of $s$ and $t$, satisfies
$\mathbf{n}_{\ell}(s,t) \leq \dst\frac{t-s}{\tau_{\ell}(a_{\ell},b_{\ell},\Delta)} +1$. Noting that ${\underset{\ell\in\{1,\ldots,M\}}\sum}\mathbf{n}_{\ell}(s,t)=k-i$,
\begin{equation}
\begin{array}{llllllllllll}
k-i & \leq  & {\underset{\ell\in\{1,\ldots,M\}}\sum}\left(\dst\frac{t-s}{\tau_{\ell}(a_{\ell},b_{\ell},\Delta)} +1\right)
\\ & \leq & M\left({\underset{\ell\in\{1,\ldots,M\}}\min}\tau_{\ell}(a_{\ell},b_{\ell},\Delta)\right)^{-1}(t-s)+ M.
\end{array}
\end{equation}
As a result, we conclude that the solutions to (\ref{eq-sys-etc-q-aux}) have a uniform semiglobal average dwell-time with $\tau(\Delta)=M^{-1}{\underset{\ell\in\{1,\ldots,M\}}\min}\tau_{\ell}(a_{\ell},b_{\ell},\Delta)$ and $n_{0}=M$ in view of Definition \ref{def-prel-adt}.

\subsection{Hybrid invariance principle}\label{subsect-proof-inv-principle}

We now apply an invariance principle for hybrid systems, namely Theorem 8.2 in \cite{Goebel-Sanfelice-Teel-book}. We introduce $U_{\aux}:\Rl{n_{q_\aux}}{}\to \Rlo$ which takes the same values as $U$ (the only difference with $U$ is its domain of definition, namely $\Rl{n_{q_\aux}}{}$ instead of \footnote{In the function $U$, the variable $p$ appears only in the form $(D^{\mathrm{T}} \otimes \mathbb{I})p$ which is replaced by $z$ in $U_{\aux}$.} $\Rl{n_q}{}$).

We deduce from (\ref{eq-proof-thm-U-flow}) and (\ref{eq-proof-thm-U-jump})
\begin{equation}
\label{eq-proof-thm-U}
\ba{rlllllll}
\langle \nabla U_\aux(q_\aux), f_\aux(q_\aux)\rangle & \leq &  u_{c}(q_\aux) &  \forall q_\aux\in C_\aux\\
U_\aux(g_\aux(q_\aux)) - U_\aux(q_\aux) & \leq&  u_{d}(q_\aux) &  \forall q_\aux\in D_\aux,
\ea
\end{equation}
where
\begin{equation}
\begin{array}{llllllllll}
u_{c}(q_\aux) & = & \left\{\begin{array}{ll}
                 -(1-\kappa)\sum_{i\in \mathcal{I}} \rho_i(v_i) & \text{when } q_\aux\in C_\aux\\
                 -\infty & \text{otherwise}
               \end{array}\right. \\
u_{d}(q_\aux) & = & \left\{\begin{array}{ll}
                 0 & \text{when } q_\aux\in D_\aux\\
                 -\infty & \text{otherwise.}
               \end{array}\right.
\end{array}\label{eq-proof-thm-uc-ud}
\end{equation}
We note that $u_c$ and $u_d$ are non-positive and that $U_\aux$ is continuous as required by Theorem 8.2 in \cite{Goebel-Sanfelice-Teel-book}. Moreover, we have shown that any maximal solution to (\ref{eq-sys-etc-aux}) is precompact. As a consequence, any maximal solution to (\ref{eq-sys-etc-aux}) approaches the largest weakly invariant subset $\mathcal{S}$ of
\begin{equation}
\begin{array}{llllll}
U_\aux^{-1}(r)\cap \mathcal{V} \cap \left[\overline{u_{c}^{-1}(0)}\cup\left(u_{d}^{-1}(0)\cap g(u_{d}^{-1}(0))\right)\right]
\end{array}
\end{equation}
where $\mathcal{V}:= \Rl{n_{q_\aux}}{}$ and $r\in U_\aux(\mathcal{V})$. Since $\overline{u_{c}^{-1}(0)}=\{q_\aux\,:\, q_\aux\in C_\aux \text{ and } v=\mathbf{0}\}$  (as $\rho_i$ is positive definite for any $i\in\mathcal{I}$, see Assumption \ref{ass-strict-passivity})  and $u_{d}^{-1}(0)=D_\aux$ in view of (\ref{eq-proof-thm-uc-ud}), the set above is
\begin{equation}
\begin{array}{c}
U_\aux^{-1}(r) \cap \left[\{q_\aux\,:\, q_\aux\in C_\aux \text{ and } v=\mathbf{0}\} \right.\\[2mm]
\left.\cup\left(D_\aux\cap g(D_\aux)\right)\right].
\end{array}\label{eq-proof-thm-set}
\end{equation}
Let $\xi\in \mathcal{S}$ and $q_\aux$ be a maximal solution such that $q_\aux(0,0)=\xi$ and $q_\aux(t,k)\in\mathcal{S}$ for any $(t,k)\in\dom q_{\aux}$, which exists as $\xi\in \mathcal{S}$ and $\mathcal{S}$ is weakly forward invariant (see Definition \ref{def-fi}).
We proceed by contradiction to show that $v(0,0)=\mathbf{0}$. Suppose $\xi\notin \overline{u_{c}^{-1}(0)}$, necessarily $\xi\in u_{d}^{-1}(0)\cap g(u_{d}^{-1}(0))= D_\aux\cap g(D_\aux)$.
The solution $q_\aux$ experiences a finite number of jumps $m\in\{1,\ldots,M\}$ until all the clocks which are equal to their lower bound are reset (and  all the variables $\hat z_\ell$ with indices $\ell$ corresponding to the clocks that were reset are updated to $z_\ell$).  After the jumps,  $q_\aux(0,m)\in C_\aux\backslash D_\aux$ in view of (\ref{eq-sys-etc-aux}). This implies that $q_\aux(0,m)\in \overline{u_{c}^{-1}(0)}$ as otherwise $q_\aux$ will no longer be in the set (\ref{eq-proof-thm-set}) (which is not possible as $q_\aux(t,k)\in \mathcal{S}$ for any $(t,k)\in\dom q_\aux$). As a consequence $v(0,m)=\mathbf{0}$ and, since  $v$ is not affected by jumps in view of (\ref{eq-sys-etc-aux}), $v(0,m)=v(0,0)=\mathrm{0}$, which contradicts the original claim that $v(0,0)\ne \mathbf{0}$. As a result $\xi\in \overline{u_{c}^{-1}(0)}$. Next we prove that $\xi\in \mathcal{S}$ and $q_\aux(0,0)=\xi$ implies $z(0,0)\in \mathcal{A}$.
In view of (\ref{eq-sys-etc-aux}), $q_\aux$ flows for at least $\varepsilon>0$ units of ordinary times from $(0,m)$ to $(\varepsilon,m)$. Consequently, for almost all $t\in[0,\varepsilon]$,
\begin{equation}
\begin{array}{lllll}
\dot{z} & = & (D^{\mathrm{T}}\otimes\mathbb{I})h(\mathbf{0})\\
\mathbf{0} & = & f(\mathbf{0},\hat u) \\
\dot{\hat{z}} &=& 0\\
\dot{\phi} & = & -\Sigma^{-1}
\left(\mathbf{1}+ \Phi(z)\right).
\end{array}\label{eq-N-sys-2nd-order-etc-hybrid.z.cont.inv.set}
\end{equation}
We remark that the vector $z$ is constant from $(0,m)$ to $(\varepsilon,m)$ as $h$ cancels at the origin in view of Section \ref{sect-pb-statement}.
Without loss of the generality, let the state $q_\aux$ stop flowing at $(\varepsilon,m)$, i.e.~ $q_\aux(\varepsilon,m)\in D_{\aux}$.  There will be again a finite number $\mu \le M$ of jumps until all the clocks which are equal to their lower bound are reset. In general, these clocks may be different from those that updated their values at the times $(0,0),\ldots,(0,m)$.  Similarly all the components of $\hat z$ corresponding to these clocks will be reset to the corresponding components of $z$.  At time $(\varepsilon,m+\mu)$, $q_\aux$ belongs to $C_\aux\backslash D_\aux$ and it starts flowing again. The clock variable $\phi_\ell$, $\ell\in\{1,\ldots,M\}$, is monotonically decreasing with a decrease rate that is bounded away from zero.  Hence, when $q_\aux$ flows, each $\phi_\ell$ decreases until eventually reaches the value $a_\ell$. There exists $\overline{M}\in\Zp$ such that, after at most $\overline{M}$ intervals during which $q_\aux$ flows, all the clock variables have undergone a reset and correspondingly all the components of $\hat z$ have been reset to the corresponding values of the components of $z$. As a result, we denote by $(\bar t, \bar k)$ the first time at which all the components of  $\hat z$ have been reset to $z$. At this time  $\hat z(\bar t, \bar k)=z(\bar t, \bar k)$. Since $\dot z=\mathbf{0}$ and $z^+=z$, then $z(\bar t, \bar k)=z(0,0)$. As a consequence, $\hat z(\bar t, \bar k)=z(0,0)$.  Let\footnote{Note that the solution may jump a finite number of times from $(\bar t,\bar k)$ to $(\bar t,\bar k')$ before flowing again, that is the reason why we consider $\bar k'$ and not $\bar k$ here.} $\bar k'\geq \bar k$ with $(\bar t,\bar k')\in\dom q_{\aux}$ be such that there exists a time $\delta>0$ such that  $q_\aux$ flows from $(\bar t, \bar k')$ to $(\bar t+\delta, \bar k')$ according to (\ref{eq-N-sys-2nd-order-etc-hybrid.z.cont.inv.set}). Since $\dot{\hat z}=\mathbf{0}$, we have that $\hat z(t,\bar k')=z(0,0)$ for all $t\in [\bar t, \bar t+\delta]$. On the other hand, the identity  $\mathbf{0} = f(\mathbf{0},\hat u)$ implies that $\hat u=\mathbf{0}$ in view of Section \ref{sect-pb-statement}. Hence $\hat u(t,\bar k')=-(D\otimes \mathbb{I})\Psi(\hat z(t,\bar k')) = \mathbf{0}$ for any $t\in[\bar t,\bar t+\delta]$, which holds during flows, entails that $\Psi(\hat z(t,\bar k'))$ belongs to the null space of $D\otimes \mathbb{I}$. Therefore $\Psi(z(0,0))$ belongs to the null space of $D\otimes \mathbb{I}$, which implies that $z(0,0)\in\mathcal{A}$ in view of item (iii) of Theorem \ref{th-etc-N}.

The arguments above show that, if $\xi\in \mathcal{S}$ and $q_\aux$ is a complete solution such that $q_\aux(0,0)=\xi$, the weak forward invariance of  $\mathcal{S}$ implies that $v(0,0)=\mathbf{0}$, $z(0,0)\in\mathcal{A}$, thus proving that
$\mathcal{S}\subseteq \mathcal{T}:=\{q_{\aux}\,:\, z\in\mathcal{A},\,\,\,v=\mathbf{0},\text{ and }\phi_{\ell}\in[a_{\ell},b_{\ell}] \text{ for all } \ell\in \{1,\ldots,M\}\}$. We conclude that any maximal solution to \eqref{eq-sys-etc-aux} approaches the largest weakly invariant set contained in $\mathcal{S}$ which is  included in $\mathcal{T}$. Consequently, any maximal solution to \eqref{eq-sys-etc} approaches the set $\{q\,:\, z\in\mathcal{A},\,\,\,v=\mathbf{0},\text{ and }\phi_{\ell}\in[a_{\ell},b_{\ell}] \text{ for all } \ell\in \{1,\ldots,M\}\}$, which concludes the proof.

\section{Conclusions}\label{sect-conclusions}

We have presented a method to design distributed controllers which ensure the coordination within a network of systems in a cyber-physical environment. Several scenarios have been investigated depending on the considered resource constraints, which we have translated as different sampling paradigms. One of the originalities of our approach is the use of auxiliary variables to define the sampling rules. This technique allows us to address a fairly general class of nonlinear networked systems, which can even be heterogeneous in the case of event-triggered and time-triggered control. The analysis is based on the hybrid formalism of \cite{Goebel-Sanfelice-Teel-book} and we have used an hybrid invariance principle to prove that the desired coordination is achieved.

A key assumption in our work is the strict passivity of the $v_i$-systems, $i\in\mathcal{I}$. This property may be ensured by an internal feedback loop in some cases. We will investigate in future work the sampling of this loop using similar techniques as those employed in this paper. 
The presented work can also serve as a basis to address other coordination problems, like when the network topology is time-varying for instance, or when the $v_i$'s have to follow a prescribed time-varying trajectories as mentioned in Remark \ref{rem-velocity-requirement}. 
Another interesting problem occurs when the reference signal for the velocities is the same for all the agents but {\em not} known to all of them. In this case, each agent should reconstruct the reference from available measurements and the problem becomes challenging even in the presence of a constant reference. This problem should be tackled relying on distributed output regulation theory for passive systems as in \emph{e.g.}, \cite{Bai2011,burger.depersis.aut13,depersis.jaya.TCNS14}.

\bibliographystyle{plain}
\bibliography{IEEEabrv,bib_global}

\end{document}

%% file: commands.tex
\newcommand{\aux}{\ensuremath{\text{aux}}}

\newcommand{\Rl}[2]{\ensuremath{\mathbb{R}^{#1}_{#2}}}   
\newcommand{\Rlp}{\ensuremath{\mathbb{R}_{>0}}}
\newcommand{\Rlo}{\ensuremath{\mathbb{R}_{\geq 0}}}
\newcommand{\Zo}{\ensuremath{\mathbb{Z}_{\geq 0}}}
\newcommand{\Zp}{\ensuremath{\mathbb{Z}_{> 0}}}

\newcommand{\eg}{{\it e.g. }}

\newcommand{\K}{\ensuremath{\mathcal{K}}}
\newcommand{\Kinf}{\ensuremath{\mathcal{K}_{\infty}}}

\newcommand{\dom}{\ensuremath{\text{dom}\,}}

\newcommand{\rge}{\ensuremath{\text{rge}\,}}
\newcommand{\diag}{\ensuremath{\text{diag}\,}}

\newcommand{\ie}{{\it i.e. }}

\newtheorem{defn}{Definition}

\newtheorem{ass}{Assumption}

\newtheorem{ex}{Example}
\newtheorem{thm}{Theorem}
\newtheorem{cor}{Corollary}

\newtheorem{rem}{Remark}

\newcommand{\dst}{\displaystyle}
\newcommand{\Ni}{\mathcal{N}_i}

\newtheorem{remark}{Remark}

\def\be{\begin{equation}}
\def\ee{\end{equation}}
\def\ba{\begin{array}}
\def\ea{\end{array}}